\documentclass[11pt]{article}
\usepackage[T1]{fontenc}
\usepackage{lmodern}
\usepackage{amsmath,amssymb,amsthm,mathtools}
\usepackage[margin=1in]{geometry}
\usepackage{microtype}
\usepackage[hidelinks]{hyperref}
\newtheorem{theorem}{Theorem}[section]
\newtheorem{proposition}[theorem]{Proposition}
\newtheorem{lemma}[theorem]{Lemma}
\newtheorem{corollary}[theorem]{Corollary}
\theoremstyle{remark}
\newtheorem{remark}[theorem]{Remark}
\theoremstyle{plain}
\DeclareMathOperator{\Ran}{Ran}

\DeclareMathOperator{\Ric}{Ric}
\DeclareMathOperator{\diver}{div}
\title{Chern scalar curvature and finite-time singularities of the conformal Chern--Calabi flow}
\author{%
Hongyu Chen\thanks{School of Mathematics, Sichuan University,
Chengdu 610064, China.
Email: \texttt{hongyu.chern@gmail.com}.}
\qquad
Li Ma\thanks{%Corresponding author.
School of Mathematics and Physics,
University of Science and Technology Beijing,
Beijing 100083, China.
Email: \texttt{nuslma@126.com}.}
\qquad
Li Sheng\thanks{School of Mathematics, Sichuan University,
Chengdu 610064, China.
Email: \texttt{lsheng@scu.edu.cn}.}
}
\date{October 2026}

\begin{document}
\maketitle

\begin{abstract}
On every unit-volume flat complex three-torus, we construct a solution of
the volume-preserving conformal Chern--Calabi flow with a finite-time
Chern Type-I singularity at any prescribed point. The full
spatial supremum has an exact asymptotic coefficient, and prescribed
rescalings converge to a complete nonflat radial shrinking solution.
The positive part of Chern scalar curvature can remain arbitrarily
small throughout the flow, as can each prescribed subcritical
spatial curvature norm. We determine the punctured terminal metric,
its one-point completion and unique nonflat tangent cone, and the
terminal energy asymptotics. We also exclude Chern Type-II under a logarithmically improved
critical weak $L^3$ bound on Chern scalar curvature, and prove a
separate scale-invariant integral criterion involving first
Chern--Ricci curvature. The constructed singular solutions satisfy
both criteria.

\smallskip

\smallskip


\end{abstract}

\tableofcontents{}

\section{Introduction}\label{sec:introduction}
  In this paper, we construct the finite time singularity of the conformal Chern-Calabi flow introduced in  \cite{MS26}
  on every flat complex three-torus, as announced in \cite{CMS26}. 
  We also exclude Chern Type-II singularities under a logarithmically improved
critical weak $L^3$ bound on the Chern scalar curvature, and prove a
separate scale-invariant integral criterion involving the first Chern--Ricci curvature.

Let $(M^3,\omega_0,J)$ be a flat complex three-torus of volume one.
For a real smooth function $f$, write $\omega_f=e^f\omega_0$ and
$\gamma_f=e^f\gamma_0$, where $\gamma_0$ is the associated real metric.
In holomorphic coordinates, write
$$
\omega_f=\sqrt{-1}\,g_{f,i\bar j}\,dz^i\wedge d\bar z^j,
\qquad (g_f^{i\bar j})=(g_{f,i\bar j})^{-1}.
$$
The first Chern--Ricci tensor and Chern scalar curvature are given by
$$
\Ric^{\rm Ch}_{i\bar j}(\omega_f)
=-\partial_i\partial_{\bar j}\log\det(g_{f,k\bar\ell}),
\qquad R_f=g_f^{i\bar j}\Ric^{\rm Ch}_{i\bar j}(\omega_f).
$$
Set $d\mu_f=\omega_f^3/3!$ and
$L_f=g_f^{i\bar j}\partial_i\partial_{\bar j}$, and let $L_f^*$
denote the formal $L^2(M,d\mu_f)$ adjoint of $L_f$.
For any smooth real $f$, define
$E(f)=\int_MR_f^2\,d\mu_f$; along a solution we write
$E(t)=E(f(t))$. We consider the volume-preserving conformal
Chern--Calabi flow
\begin{equation}\label{eq:flow}
f_t=6L_f^*R_f-R_f^2+E(t).
\end{equation}
The Laplacian is nonpositive at a local maximum. The first variation
in Section~\ref{sec:equation} gives
\begin{equation}\label{eq:dissipation}
\int_Md\mu_f=1,\qquad E'(t)=-\int_Mf_t^2\,d\mu_f.
\end{equation}

Spatial curvature norms use the evolving metric $\omega_f(t)$ and
volume form $d\mu_f(t)$. Pointwise norms of Chern curvature tensors
are the Hermitian Hilbert--Schmidt norms induced by $g_{f,i\bar j}$;
real covector norms use $\gamma_f(t)$. At each fixed time,
$$
\|R_f(t)\|_\infty
=\|R_f(t)\|_{L^\infty(M,d\mu_f(t))}
=\max_{x\in M}|R_f(t,x)|.
$$
For tensors, the spatial $L^\infty$ norm is the maximum of the
pointwise norm induced by $\omega_f(t)$. Time suprema are written
explicitly.

The conformal transformation of Chern scalar curvature underlies the
Chern--Yamabe problem of Angella, Calamai and Spotti~\cite{ACS17}.
Equation~\eqref{eq:flow}
is instead the constrained negative gradient flow of the squared
Chern scalar curvature in complex dimension three.
Sheng and Ma~\cite{MS26} prove global existence on compact complex
surfaces with Gauduchon backgrounds and announce higher-dimensional
finite-time blow-up. The conformal weight of the squared energy is
constant in complex dimension two and equals $e^f$ here.
The precise coefficient convention in their general-dimensional
formula is discussed in Section~\ref{sec:equation}.
The potential flows of Bedulli--Vezzoni~\cite{BV20} and
Shen~\cite{Shen22} evolve different Hermitian data.

Curvature bounds also govern the regularity of other Calabi-type flows.
For the K\"ahler Calabi flow, Chen and He~\cite{CH08} prove that a
uniform Ricci curvature bound gives extension across a finite time.
Kawamura~\cite{Kawamura26} obtains a priori estimates for the Calabi
flow in a quasi-K\"ahler setting under a uniform Chern scalar
curvature bound.

Chen, Ma and Sheng~\cite[Theorems~1.1--1.3]{CMS26} prove continuation
if $\int_0^T\|R_f(t)\|_{L^\infty(M)}^2\,dt<\infty$, and global convergence
under a uniform evolving $L^p$ curvature bound for one $p>3$.
Their Introduction announces the finite-time nonextension
construction treated here. We prove the construction on each fixed
flat torus, identify its scalar curvature rate and limiting geometry,
and examine curvature conditions that permit this behavior.

For fourth-order curvature flows, Streets~\cite{Streets13long}
obtains smoothing estimates under a curvature bound. For the $L^2$
curvature flow, he proves that the curvature supremum is unbounded
at a finite-time singularity. He also constructs finite-time
singularities in dimensions at least five, including examples in
dimensions at least six that collapse at the curvature
scale~\cite{Streets13collapse}. These results provide a comparison
for the singularity formation and terminal geometry studied here.

Under fourth-order geometric scaling, the metric is multiplied by
$K$, time by $K^2$, and Chern scalar curvature by $K^{-1}$.
We classify finite-time singularities by the growth of Chern scalar
curvature. A singularity at a finite maximal time $T$ is called
\emph{Chern Type-I} if
\begin{equation}\label{eq:typeI}
\sup_{0\le t<T}\sqrt{T-t}\,\max_{x\in M}|R_f(t,x)|<\infty,
\end{equation}
and \emph{Chern Type-II} if the corresponding terminal limsup is infinite.
The full Chern tensor is denoted by $\mathcal R^{\rm Ch}$. The notation
$\nabla=\nabla^{\rm Ch}$ always means the evolving Chern connection;
$\nabla^0$ denotes the fixed flat background Chern connection.

\begin{theorem}\label{thm:construction}
Let $(M^3,\omega_0,J)$ be any flat complex three-torus of volume one
and let $p\in M$. There are smooth unit-volume initial data whose
maximal solution of \eqref{eq:flow} has finite maximal time
$\Theta>0$, with
$$
0<\lim_{t\uparrow\Theta}\sqrt{\Theta-t}\,\max_{x\in M}|R_f(t,x)|<\infty.
$$
The only singular point is $p$, where $e^{f(t,p)}\to0$ and
$R_f(t,p)\to-\infty$. On $M\setminus\{p\}$ the metrics converge
smoothly to a positive metric $\gamma_\Theta$.
\end{theorem}

For Ricci flow, Enders, M\"uller and Topping~\cite{EMT11} show
that pointed blow-up limits at singular points of Type-I flows are
nontrivial gradient shrinking solitons.
The shrinking model is defined in Section~\ref{sec:shrinking}.
Proposition~\ref{prop:blowup-geometry} identifies the exact curvature
coefficient and the smooth rescaled limit.
Corollary~\ref{cor:terminal} determines the terminal completion,
its unique tangent cone and Gromov--Hausdorff convergence.

\begin{corollary}\label{cor:smallness}
Let $(M^3,\omega_0,J)$ be any flat complex three-torus of volume one
and let $p\in M$. For every $1\le p_0<3$ and
$\varepsilon,\tau_*>0$, there are smooth unit-volume initial data
whose maximal solution of \eqref{eq:flow} has finite maximal time
$0<\Theta<\tau_*$, satisfies all conclusions of
Theorem~\ref{thm:construction}, and obeys
\begin{equation}\label{eq:smallpositive}
R_f(t,x)<\varepsilon\qquad(0\le t<\Theta,\ x\in M).
\end{equation}
For every $0\le t<\Theta$, one also has
\begin{equation}\label{eq:smallsubcritical}
\|R_f(t)\|_{L^{p_0}(M,d\mu_f(t))}
+\|\Ric^{\rm Ch}(\omega_f(t))\|_{L^{p_0}(M,d\mu_f(t))}
+\|\mathcal R^{\rm Ch}(t)\|_{L^{p_0}(M,d\mu_f(t))}
<\varepsilon.
\end{equation}
\end{corollary}

These solutions have $R_f(t,p)\to-\infty$ despite the upper bound
\eqref{eq:smallpositive}. Thus a scalar curvature upper bound,
even arbitrarily small throughout the flow, does not give continuation.
The subcritical norms in \eqref{eq:smallsubcritical} contrast with
the supercritical continuation theorem in~\cite{CMS26}.
The critical norm in this construction has a logarithmic divergence,
computed in Proposition~\ref{prop:critical}.

Our scalar curvature criterion uses the weak Lorentz quasi-norm
$$
\|v\|_{L^{3,\infty}(M,d\mu)}^3
=\sup_{\lambda>0}\lambda^3\mu\{x\in M:|v(x)|>\lambda\}.
$$
Distribution functions and critical weak Lorentz spaces also occur
in the geometric analysis of Ricci shrinkers. Yudowitz uses this
framework in the analysis of bubble-tree convergence and Morse
index, with critical weights in $L^{n/2,\infty}$; see
\cite[Definition~4.1 and Proposition~4.6]{Yudowitz25}.
Here the real dimension is six, and the critical curvature exponent
is three.

\begin{theorem}\label{thm:scalar-criterion}
Let $f$ be a maximal smooth unit-volume solution of \eqref{eq:flow}
with finite maximal time $T$. If, for some $0\le t_0<T$,
\begin{equation}\label{eq:scalar-log}
\sup_{t_0<t<T}
\left\|R_f(t)
\left[\log\!\left(e+\sqrt{T-t}\,|R_f(t)|\right)\right]^{1/3}
\right\|_{L^{3,\infty}(M,d\mu_f(t))}<\infty,
\end{equation}
then the solution is Chern Type-I.
\end{theorem}

\begin{remark}\label{rem:scalar-gauge}
More generally, let $h:[0,\infty)\to[1,\infty)$ be nondecreasing
and locally bounded, with $h(s)\to\infty$ as $s\to\infty$.
The same conclusion holds if
\begin{equation}\label{eq:scalar-gauge}
\sup_{t_0<t<T}\sup_{\lambda>0}
\lambda^3h\!\left(\sqrt{T-t}\,\lambda\right)
\mu_{f(t)}\{|R_f(t)|>\lambda\}<\infty.
\end{equation}
Theorem~\ref{thm:scalar-criterion} is the case $h(s)=\log(e+s)$.
One may also take $h(s)=[\log(e+s)]^\varepsilon$ for any
$\varepsilon>0$, or $h(s)=\log\log(e^e+s)$.
\end{remark}

Under fourth-order scaling, $|R_f|^3d\mu_f$ and
$\sqrt{T-t}|R_f|$ are invariant. The logarithm in
\eqref{eq:scalar-log} is comparable to a constant below the Type-I
curvature scale and increases above it. Proposition~\ref{prop:compatibility}
shows that the constructed singular solutions satisfy
\eqref{eq:scalar-gauge} for every $h$ in
Remark~\ref{rem:scalar-gauge}, although their strong $L^3$ curvature
integrals diverge. The proof of the scalar criterion combines a
static conformal gradient estimate with local regularity on the
past of a first curvature maximum.

We also retain a scale-invariant terminal integral criterion for
the first Chern--Ricci curvature.

\begin{theorem}\label{thm:ricci-criterion}
Let $f$ be a maximal smooth unit-volume solution of \eqref{eq:flow}
with finite maximal time $T$, and fix $0<s_0<T/2$. If
\begin{equation}\label{eq:ricci-budget}
\sup_{0<s<s_0}\int_{T-2s}^{T-s}
\|\Ric^{\rm Ch}(\omega_f(t))\|_{L^9(M,d\mu_f(t))}^{3}\,dt<\infty,
\end{equation}
then the solution is Chern Type-I and
$$
\sup_{0\le t<T}\sqrt{T-t}\,
\|\Ric^{\rm Ch}(\omega_f(t))\|_{L^\infty(M)}<\infty.
$$
\end{theorem}

Proposition~\ref{prop:compatibility} proves that the constructed
singular solutions have a positive finite limiting budget in
\eqref{eq:ricci-budget}. Under fourth-order scaling the spatial
$L^9$ curvature norm is multiplied by $K^{-2/3}$, so the displayed
time integral is invariant. The proof first controls the conformal
gradient by this curvature norm; local parabolic estimates then
supply time propagation. Whether the terminal dyadic integral of
$\|R_f\|_\infty^2$ or $\|\Ric^{\rm Ch}\|_\infty^2$ alone implies
\eqref{eq:typeI} remains unresolved.

The construction combines a smooth radial shrinking profile with
a finite-dimensional selection on the fixed torus.
Weighted estimates in the rescaled core are coupled to estimates
on physical annuli and the exterior. A Newton potential permits
the negative Laplacian of the exterior conformal factor to be
arbitrarily small. The global curvature rate follows from a uniform
tail estimate, and the original time normalization is handled
directly in affine cylinders.

Here is the plan of the paper.
Section~\ref{sec:equation} derives the equation and gives the local
analytic estimates. Section~\ref{sec:criteria} proves the scalar
and Chern--Ricci curvature criteria. Sections~\ref{sec:shrinking} and
\ref{sec:localization} construct and localize the shrinking profile.
Section~\ref{sec:smallness} proves the full-time smallness assertions,
and Section~\ref{sec:constructed-geometry} determines the rate,
the limiting geometry and terminal energy.
Appendices~\ref{sec:profile} and \ref{sec:linear} contain the profile
matching and the weighted linear analysis.

\section{The equation and local estimates}\label{sec:equation}

Write $d\mu_0=\omega_0^3/3!$ and, in holomorphic coordinates,
$L_0=g_0^{i\bar j}\partial_i\partial_{\bar j}$. Then
$$
L_f=e^{-f}L_0,\qquad R_f=-3e^{-f}L_0f,\qquad
d\mu_f=e^{3f}d\mu_0.
$$
For a real variation $\varphi$,
$$
\left.\partial_\epsilon R_{f+\epsilon\varphi}\right|_{\epsilon=0}
=-\varphi R_f-3L_f\varphi,\qquad
\left.\partial_\epsilon d\mu_{f+\epsilon\varphi}\right|_{\epsilon=0}
=3\varphi\,d\mu_f.
$$
Consequently
\begin{equation}\label{eq:variation}
DE(f)[\varphi]=
\int_M(R_f^2-6L_f^*R_f)\varphi\,d\mu_f.
\end{equation}
The tangent space to unit volume consists of $\varphi$ with
$\int_M\varphi\,d\mu_f=0$. Since
$\int_ML_f^*R_f\,d\mu_f=0$, projection of the negative gradient
in \eqref{eq:variation} gives \eqref{eq:flow}.
Substitution of $\varphi=f_t$ proves \eqref{eq:dissipation}.
The leading adjoint coefficient in formula (4) of~\cite{MS26} is
$2$, whereas the first variation in complex dimension $n$ gives
$2n$. In dimension two the reaction term vanishes and the flows
are related by a constant time change. In higher dimensions their
relative coefficients differ; we use \eqref{eq:flow}.

Two integrations by parts in the fixed flat measure give
\begin{equation}\label{eq:adjoint}
L_f^*\varphi=e^{-3f}L_0(e^{2f}\varphi).
\end{equation}
Choose affine real coordinates $X=(X^1,\ldots,X^6)$ with
$\gamma_0=\tfrac12g_{\mathrm E}$, where
$g_{\mathrm E}=\sum_{a=1}^6(dX^a)^2$.
We use $\Delta_X=\sum_{a=1}^6\partial_{X^a}^2$ and the sign
$\Delta_{\gamma_0}=\operatorname{div}_{\gamma_0}
\operatorname{grad}_{\gamma_0}$ for the real Laplace--Beltrami operator.
Thus
$$
L_0=\Delta_X=\tfrac12\Delta_{\gamma_0},\qquad d\mu_0=dX/8.
$$
The affine derivatives used below are components of $\nabla^0$,
and $\nabla^{0,j}=(\nabla^0)^j$ denotes ordered derivatives.
In affine formulae, an unsubscripted $\Delta$ is the Euclidean
Laplacian in the current variables; subscripts $X$, $y$ and $Y$
indicate those variables when needed. Euclidean metrics in these
variables are denoted by $g_{\mathrm E}$.
Gradients, divergences and coefficient contractions in affine
formulae use $g_{\mathrm E}$, as do the norms and inner products
with subscript $E$. Intrinsic norms carry their displayed metric.
In particular,
$|df|_{\gamma_f}^2=2e^{-f}|\nabla^0f|_E^2$.
This distinction preserves the factors in the complex Laplacian.

The unnormalized local equation is $\widetilde f_s=\mathcal N(\widetilde f)$,
where
\begin{equation}\label{eq:N}
\begin{split}
\mathcal N(f)
&=-18e^{-3f}\Delta_X(e^f\Delta_Xf)-9e^{-2f}(\Delta_Xf)^2\\
&=-18e^{-2f}\left(\Delta_X^2f
+2\langle\nabla^0f,\nabla^0\Delta_Xf\rangle_E
+|\nabla^0f|_E^2\Delta_Xf+\tfrac32(\Delta_Xf)^2\right).
\end{split}
\end{equation}
The same formula in affine model variables will also be denoted by
$\mathcal N$. It has the constant-shift law
$\mathcal N(f+c)=e^{-2c}\mathcal N(f)$.

\begin{proposition}\label{prop:existence}
Every smooth real $f_0$ with $\int_Me^{3f_0}d\mu_0=1$ has a
unique maximal smooth solution of \eqref{eq:flow} on
$[0,T_{\max})$, where $T_{\max}>0$.
\end{proposition}

\begin{proof}
The unnormalized principal part is $-18e^{-2f}L_0^2f$,
or $-\tfrac92e^{-2f}\Delta_{\gamma_0}^2f$.
Its symbol is the negative square of a positive quadratic form.
The other terms have order at most three and smooth coefficients.
The short-time theorem of Mantegazza--Martinazzi~\cite[Theorem~1.1]{MM12}
applies on the closed manifold.

For its solution put $V(s)=\int_Me^{3\widetilde f(s)}d\mu_0$.
Integration of the equation yields $V'=-3E(\widetilde f)$.
On the interval where $V>0$, define
\begin{equation}\label{eq:normalize}
t_s=V^{-2/3},\qquad
f(t(s))=\widetilde f(s)-\tfrac13\log V(s),\qquad t(0)=0.
\end{equation}
Curvature and volume scaling give
$E(f)=V^{-1/3}E(\widetilde f)$. Differentiation of
\eqref{eq:normalize} therefore gives
$f_t=V^{2/3}\widetilde f_s+E(f)=\mathcal N(f)+E(f)$.
Conversely, a normalized solution is transformed back by subtracting
$\int_0^tE(r)\,dr$ and setting
$ds/dt=\exp(-2\int_0^tE(r)\,dr)$.
Uniqueness for the unnormalized equation gives uniqueness after
normalization. The union of compatible local solutions gives the
maximal solution.
\end{proof}

\begin{lemma}\label{lem:Sobolev-continuation}
On a fixed flat complex three-torus, a smooth solution of
$\widetilde f_s=\mathcal N(\widetilde f)$ extends across a finite
endpoint if its $H^{14}$ norm in the fixed background metric remains
bounded there. The local existence interval and continuous
dependence in $H^{14}$ can be chosen uniformly on bounded sets of
smooth initial data in this norm.
\end{lemma}

\begin{proof}
For this proof write $f$ for the unnormalized unknown. Sobolev
embedding in real dimension six bounds its $C^{10}$ norm by its
$H^{14}$ norm. The differential expression has the divergence identity
\begin{align*}
\mathcal N(f)={}&-18\Delta(e^{-2f}\Delta f)-108e^{-2f}\nabla^0 f\cdot\nabla^0\Delta f\\
&-63e^{-2f}(\Delta f)^2+54e^{-2f}|\nabla^0 f|^2\Delta f.
\end{align*}
If $\|f\|_{H^{14}}\le2K$, the coefficient $e^{-2f}$ has a positive lower bound and bounded $C^{10}$ norm. For every integer $j\ge14$, the commutator in the principal energy satisfies
$$
\|[\nabla^{0,j},e^{-2f}]\Delta f\|_2
\le C_j(K)(\|f\|_{H^{j+1}}+\|f\|_{H^j}).
$$
The other three terms have $H^{j-2}$ norm at most $C_j(K)(1+\|f\|_{H^{j+1}})$. Integration by parts twice, interpolation and Young's inequality therefore give
$$
\frac d{ds}\|f\|_{H^j}^2+c(K)\|f\|_{H^{j+2}}^2
\le C_j(K)(1+\|f\|_{H^j}^2).
$$
Here and below constants may depend on the fixed background and fixed
parameters, but never on a preterminal endpoint. For $j=14$ the
right side is bounded in terms of $K$, supplying a uniform positive
local time. Fourier Galerkin projection commutes with derivatives
and preserves this estimate. Compactness, the bound on $f_s$ in
$L^2H^{12}$, and the interpolation trace $C H^{14}$ give a solution.

To verify the continuous dependence used in the degree argument,
let $f$ and $g$ be two such solutions with $H^{14}$ norms at most
$2K$. The coefficient difference in the leading divergence term
satisfies the product estimate
$$
\|(e^{-2f}-e^{-2g})\Delta g\|_{H^{14}}
\le C(K)\|f-g\|_{H^{14}}(1+\|g\|_{H^{16}}).
$$
It follows by placing either the coefficient difference in
$H^{14}$ and $\Delta g$ in $L^\infty$, or the coefficient
difference in $L^\infty$ and $\Delta g$ in $H^{14}$.
The difference of each lower term has the same estimate in
$H^{12}$, with an additional $C(K)\|f-g\|_{H^{15}}$.
Pairing the difference equation with $f-g$ at order fourteen,
integrating the principal divergence term twice, and absorbing the
$H^{15}$ term by interpolation gives
\begin{align*}
&\frac d{ds}\|f-g\|_{H^{14}}^2
+c(K)\|f-g\|_{H^{16}}^2\\
&\qquad\le C(K)
\left(1+\|f\|_{H^{16}}^2+\|g\|_{H^{16}}^2\right)
\|f-g\|_{H^{14}}^2.
\end{align*}
The coefficient is integrable on the common local interval, by the
preceding energy estimate. Gronwall's inequality therefore gives
uniqueness and continuous dependence in $C H^{14}\cap L^2H^{16}$.
The difference equation gives continuous dependence of the time
derivative in $L^2H^{12}$ as well. The estimates for $j>14$ give
smooth persistence on this same interval for smooth data.
Restarting at times tending to a finite endpoint proves the assertion.

\end{proof}

\begin{lemma}\label{lem:smoothing}
Let $B_2\subset\mathbb R^6$ and let $v$ be a smooth real solution
on $B_2\times(-1,0]$ of
\begin{equation}\label{eq:forced}
v_s=-18e^{-3v}\Delta(e^{v}\Delta v)
-9e^{-2v}(\Delta v)^2+\zeta(s),
\end{equation}
where $\zeta$ depends only on time. Suppose
$|v|+|R_v|+|\zeta|\le H$ on this cylinder, with
$R_v=-3e^{-v}\Delta v$. For every integer $j\ge0$,
$$
\sup_{B_1\times[-1/2,0]}
\left(|\nabla^{0,j}v|+|\nabla^{0,j}v_s|\right)\le C(H,j).
$$
The same conclusion
holds on any fixed pair of nested backward cylinders, with constants
also depending on their spatial and temporal margins.
\end{lemma}

\begin{proof}
The identity $\Delta v=-e^vR_v/3$ bounds the Laplacian.
The interior Poisson gradient estimate on each time slice then
bounds $\nabla^0v$ on $B_{7/4}$ in terms of $H$.
We work on this smaller outer ball. Expansion of
\eqref{eq:forced} gives
\begin{equation}\label{eq:linearform}
\begin{split}
v_s+18e^{-2v}\Delta^2v
+36e^{-2v}\nabla^0 v\cdot \nabla^0\Delta v
={}&-27e^{-2v}(\Delta v)^2\\
&-18e^{-2v}|\nabla^0 v|^2\Delta v+\zeta(s).
\end{split}
\end{equation}
Its fourth-order symbol is $18e^{-2v}|\xi|^4$, bounded above
and below by positive constants depending only on $H$. The
fourth-order coefficients have a uniform spatial Lipschitz bound,
the third-order coefficients are bounded, and the right-hand side is
bounded.

The non-divergence estimate of Dong and Kim~\cite[Theorem~2.3 and
Remark~2.4]{DK10} applies to scalar fourth-order operators with this
positive symbol, bounded lower-order coefficients, and sufficiently
small spatial mean oscillation of their leading coefficients on small
cylinders. The spatial Lipschitz bound supplies the last condition
at a uniform radius, independently of time derivatives. Localizing
their estimate yields, for every $1<P<\infty$, a bound on an interior
backward cylinder for $v_s$ and all $\nabla^{0,j}v$, $j\le4$, in $L^P$,
in terms of the $L^P$ norms of $v$ and the right-hand side of
\eqref{eq:linearform}. We give the localization details needed here.

Multiply $v$ by a smooth cutoff that vanishes near the spatial
boundary and near the bottom time. Extend the coefficients outside
its support while retaining the uniform spatial Lipschitz bound and
ellipticity. The cutoff commutators have spatial order at most three;
the time commutator has order zero. For smooth compactly supported
spatial functions $z$, the interpolation inequalities
$$
\|\nabla^{0,j}z\|_{L^P}\le\varepsilon\|\nabla^{0,4}z\|_{L^P}
+C(P,j)\varepsilon^{-j/(4-j)}\|z\|_{L^P},
\qquad 1\le j\le3,
$$
absorb the third-order terms after integration in time and Young's
inequality. A positive
zeroth-order parameter in the cited global estimate contributes only
a fixed multiple of $\|v\|_{L^P}$. More explicitly, there is
$N=N(H,P)>0$, independent of $\varepsilon$, such that for each
$0<\varepsilon<1$ there is $C_\varepsilon=C(H,P,\varepsilon)>0$
with the following estimate on nested cylinders of spatial gap
$R-r$ and temporal gap $t_2-t_1$:
\begin{align*}
&\|v_s\|_{L^P(B_r\times(t_2,0))}
+\sum_{j=0}^4\|\nabla^{0,j}v\|_{L^P(B_r\times(t_2,0))}\\
&\quad\le\varepsilon\left(
\|v_s\|_{L^P(B_R\times(t_1,0))}
+\sum_{j=0}^4\|\nabla^{0,j}v\|_{L^P(B_R\times(t_1,0))}\right)\\
&\qquad\quad+C_\varepsilon
\left(1+(R-r)^{-N}+(t_2-t_1)^{-N}\right).
\end{align*}
Here $1<r<R<7/4$ and $-1<t_1<t_2<-1/2$. The bounded lower-order
right-hand side has been included in the last constant.
Choose radii increasing to a number less than $7/4$, and bottom times
decreasing to a number greater than $-1$, with successive gaps
proportional to $2^{-j}$. Choosing $\varepsilon$ so that
$\varepsilon2^N<1$ makes iteration summable. Its final error is
$\varepsilon^j$ times a larger-cylinder high norm. That norm is
finite for each smooth solution on the fixed compact outer cylinder,
so the error tends to zero. This establishes a uniform interior
$W_P^{1,4}$ bound without any assumed high norm at the bottom time.

Choose $P>10$. Parabolic Sobolev embedding, with spatial order
four and temporal order one, gives spatial $C^{3,\alpha}$ bounds
on a smaller cylinder for any
$0<\alpha<1-10/P$, together with the corresponding parabolic
H\"older bounds. The coefficients and the right-hand side of
\eqref{eq:linearform} now have uniform spatial $C^\alpha$ norms.
The interior
estimate of Dong and Zhang~\cite[Theorem~2.2(I)]{DZ15}, which permits
measurable time dependence and assumes uniform spatial
$C^\alpha$ coefficients and forcing, gives uniform bounds for
$\nabla^{0,4}v$ in $C^{\alpha/4,\alpha}$ and for $v_s$ in spatial
$C^\alpha$. In particular these are bounds up to the upper time
slice, obtained by taking its limit from below.

For the iteration, keep only $18e^{-2v}\Delta^2$ on the left
of the equation and put the third-order expression on the right.
Suppose the spatial derivatives through order $j+3$ have uniform
$C^\alpha$ bounds on a cylinder. Differentiate the equation $j$
times in space. Apart from
$18e^{-2v}\Delta^2\nabla^{0,j}v$, each resulting term contains derivatives
of $v$ of order at most $j+3$: in a commutator with the leading
coefficient, at least one derivative falls on $e^{-2v}$. Thus the
same estimate controls $\nabla^{0,j+4}v$. Induction on nested cylinders supplies every
spatial derivative. The equation then supplies every $\nabla^{0,j}v_s$.
Spatial derivatives of $\zeta$ vanish, and no time derivative of
$\zeta$ is used in these estimates.
For each fixed requested order the margins can be divided into
finitely many positive pieces, ending on
$B_1\times[-1/2,0]$.
\end{proof}

\section{Curvature criteria}\label{sec:criteria}

\subsection{A Chern--Ricci integral bound}

We first prove a weighted estimate on a planar domain. Its role is
to exclude a possible zero set in the conformal rescaling below.
All derivatives and measures in the planar argument are Euclidean.

\begin{lemma}\label{lem:planar-cap}
There is no sequence of smooth positive functions $u_i$ on a
neighborhood of $\overline{B_1}\subset\mathbb C$ such that
$$
u_i\longrightarrow |z|\quad\hbox{uniformly on }\overline{B_1},
\qquad
\int_{B_1}u_i^{12}|\Delta\log u_i|^9\,dA\longrightarrow0.
$$
\end{lemma}

\begin{proof}
We first consider a fixed smooth positive function $u$ on a bounded
planar domain $D$ and a regular value $s<\min_{\partial D}u$ with
$\{u<s\}\ne\varnothing$. Suppose initially that the critical points
below this level are nondegenerate. For $\min_Du<\sigma\le s$, set
$$
H(\sigma)=|\{u<\sigma\}|,\qquad
q(\sigma)=\frac1{2\pi}\int_{\{u<\sigma\}}\Delta\log u\,dA.
$$
The divergence and coarea formulas give, at every regular level,
\begin{align*}
2\pi\sigma q(\sigma)&=\int_{u=\sigma}|Du|\,d\ell,
& H'(\sigma)&=\int_{u=\sigma}|Du|^{-1}\,d\ell,\\
2\pi q'(\sigma)&=\int_{u=\sigma}\frac{\Delta\log u}{|Du|}\,d\ell.
\end{align*}
All components and holes of the sublevel set are included.
Cauchy--Schwarz and the planar isoperimetric inequality imply
$$
2\pi\sigma q H'\ge |\partial\{u<\sigma\}|^2\ge4\pi H.
$$
H\"older's inequality on the same complete level set gives
$$
\int_{u=\sigma}\frac{\sigma^{12}|\Delta\log u|^9}{|Du|}\,d\ell
\ge (2\pi)^9\sigma^{12}|q'(\sigma)|^9H'(\sigma)^{-8}.
$$
Use the strictly increasing area coordinate
$\xi=\tfrac12\log(H(\sigma)/\pi)$ and regard $\sigma$ and $q$ as
functions of $\xi$. Since $\sigma_\xi H'=2H$, the preceding
inequalities become
\begin{equation}\label{eq:area-flux}
\frac{d\log\sigma}{d\xi}\le q,\qquad
\int_{u<s}u^{12}|\Delta\log u|^9\,dA
\ge2\pi\int_{-\infty}^{\xi_s}
\sigma^{12}e^{-16\xi}|q_\xi|^9\,d\xi.
\end{equation}
These formulas hold across the finitely many critical levels.
Indeed the coarea expressions are locally integrable there, and
$$
|q_\xi|\le\frac H\pi\sup_{\{u\le s\}}|\Delta\log u|,
\qquad 0\le(\log\sigma)_\xi\le q
$$
give local Lipschitz bounds in $\xi$. At the minimum end,
$H\to0$, $q\to0$ and $\log\sigma\to\log\min_Du$.

Suppose that $0<q(s)<4/3$. Put
$w(\xi)=\log\sigma(\xi)-4\xi/3$, and take the last interval ending
at $\xi_s$ on which $q<4/3$. On that interval
$-w_\xi\ge4/3-q>0$, so $w$ is a valid decreasing coordinate.
Changing variables in \eqref{eq:area-flux} yields
\begin{align*}
\int e^{12w}|q_\xi|^9\,d\xi
&\ge\int e^{12w}(4/3-q)^8|q_w|^9\,dw\\
&=\left(\frac9{17}\right)^9
\int e^{12w}\left|\frac d{dw}(4/3-q)^{17/9}\right|^9\,dw.
\end{align*}
The integrals in $w$ run from $w(\xi_s)$ to the value at the left
endpoint. That endpoint is either a last crossing $q=4/3$, or
$\xi=-\infty$, where $q\to0$ and $w\to+\infty$.
Weighted H\"older, using
$$
\int_{w(\xi_s)}^\infty e^{-3w/2}\,dw
=\tfrac23e^{-3w(\xi_s)/2},
$$
therefore proves
\begin{align*}
\int_{u<s}u^{12}|\Delta\log u|^9\,dA
\ge{}&2\pi\left(\frac9{17}\right)^9\left(\frac32\right)^8
s^{12}(H(s)/\pi)^{-8}\\
&\cdot\min\left\{(4/3-q(s))^{17},
\left[(4/3)^{17/9}-(4/3-q(s))^{17/9}\right]^9\right\}.
\end{align*}
For an infinite endpoint this follows first on compact intervals
and then by the stated endpoint limits. Oscillations of the flux
inside the interval do not affect the lower bound by its endpoint
variation.

To remove the nondegeneracy assumption, perturb $u$ by arbitrarily
small generic linear functions. Sard's theorem applied to $Du$
makes all critical points in the relevant compact sublevel set
nondegenerate. Positivity, the strict boundary inequality and the
regular outer level persist. Their sublevel areas and fluxes
converge to those of $u$. For this fixed function the positive
minimum and $C^2$ convergence also give convergence of the weighted
energies, including on the moving sublevel domains. Thus the lower
bound holds for every smooth positive $u$.

Now suppose a sequence in the lemma exists. On a fixed annulus
around $|z|=1/2$, uniform positivity converts its weighted bound
into $\Delta\log u_i\to0$ in ordinary $L^9$.
Uniform logarithmic convergence and interior Poisson estimates
then give $C^1$ convergence near that circle. The entire sublevel
set at $s=1/2$ has a single radial graph as boundary: uniform
convergence confines its boundary to this annulus, where the radial
derivative is positive. Consequently
$H_i(s)\to\pi s^2$ and $q_i(s)\to1$. The displayed lower bound has
a strictly positive lower limit, contradicting the assumed energy
convergence.
\end{proof}

\begin{lemma}\label{lem:ricci-gradient}
There is $C=C(M,\omega_0)$ such that every smooth real $f$ satisfies
\begin{equation}\label{eq:ricci-gradient}
\|df\|_{L^\infty(\gamma_f)}^2
\le C\|\Ric^{\rm Ch}(\omega_f)\|_{L^9(M,d\mu_f)}^{3/2}.
\end{equation}
\end{lemma}

\begin{proof}
Adding a constant to $f$ multiplies both sides by the same factor,
so normalize $\int_Me^{3f}d\mu_0=1$. In the affine convention of
Section~\ref{sec:equation}, let $u=e^{-f/2}$. Direct calculation gives
$$
|df|_{\gamma_f}^2=8|\nabla^0u|_E^2,\qquad
\Ric^{\rm Ch}_{i\bar j}=6\partial_i\partial_{\bar j}\log u.
$$
Thus the ninth power of the curvature norm is a fixed positive
multiple of
$\int_Mu^{12}|\partial\bar\partial\log u|_E^9d\mu_0$.
If the estimate fails, choose normalized $f_i$ and $u_i=e^{-f_i/2}$
with nonzero curvature such that
$$
\frac{\|\nabla^0u_i\|_{L^\infty_E}^2}
{\|\Ric^{\rm Ch}(\omega_{f_i})\|_{L^9(d\mu_{f_i})}^{3/2}}
\longrightarrow\infty.
$$
Zero curvature makes $f_i$ harmonic on the compact torus and hence
constant, so those factors already satisfy the estimate.

First suppose that the gradient norms are bounded. The curvature
norms then tend to zero. Unit volume gives
$\min u_i\le1\le\max u_i$, and the gradient bound gives
$\max u_i\le1+C_M\|\nabla^0u_i\|_{L^\infty_E}$.
A subsequence converges uniformly to a nonnegative Lipschitz function
$u$ with positive maximum. On compact subsets of $\{u>0\}$,
the weighted bound controls $\partial\bar\partial\log u_i$ in
ordinary $L^9$. Its trace and the uniform amplitudes give local
$W^{2,9}$ convergence, so $\log u$ is pluriharmonic there.
In particular $u$ is plurisubharmonic on its positive set.
For every $\delta>0$, the function $\max\{u,\delta\}$ is globally
plurisubharmonic: it is constant near zeros, and near its positive
transition level it is the maximum of two plurisubharmonic functions.
Letting $\delta\downarrow0$ shows that $u$ is plurisubharmonic on
the compact torus, hence constant. Its positive maximum and volume
normalization give $u=1$.
The factors $u_i$ and their inverses are now uniformly bounded.
The fixed-background elliptic estimate for $\log u_i$ minus its mean
therefore gives
$$
\|\nabla^0u_i\|_{L^\infty_E}
\le C\|\Ric^{\rm Ch}(\omega_{f_i})\|_{L^9(d\mu_{f_i})}.
$$
The contrary quotient is bounded by a constant times the square
root of this curvature norm and tends to zero, a contradiction.

We may therefore assume
$\|\nabla^0u_i\|_{L^\infty_E}\to\infty$.
Choose $x_i$ attaining this norm and define on the periodic affine
cover
$$
w_i(y)=\frac1{u_i(x_i)}
u_i\left(x_i+
\frac{u_i(x_i)}{\|\nabla^0u_i\|_{L^\infty_E}}y\right).
$$
Then $w_i(0)=1$, $|\nabla^0w_i|_E\le1$ and $|\nabla^0w_i(0)|_E=1$.
The volume anchor gives
$u_i(x_i)/\|\nabla^0u_i\|_{L^\infty_E}\le C_M+o(1)$.
Every fixed rescaled ball therefore meets a bounded number of
rescaled fundamental cells. Exact change of variables yields
$$
\int_{B_R}w_i^{12}|\partial\bar\partial\log w_i|_E^9dy
\le C_R
\frac{\|\Ric^{\rm Ch}(\omega_{f_i})\|_{L^9(d\mu_{f_i})}^{9}}
{\|\nabla^0u_i\|_{L^\infty_E}^{12}}
\longrightarrow0.
$$
No expansion of the rescaled periods is required.

Take a locally uniform Lipschitz limit $w\ge0$, with $w(0)=1$.
On the positive component containing zero, the preceding energy
bound and interior Poisson estimates give local $W^{2,9}$
convergence of $\log w_i$, hence $C^1$ convergence in real
dimension six. Thus $\log w$ is pluriharmonic and $|\nabla^0w(0)|_E=1$.
Locally $w=|g|$ for a nonvanishing holomorphic function, and
$$
|\nabla^0w|_E^2=\sum_{j=1}^3|g_j|^2\le1.
$$
This subharmonic function attains its maximum at zero. Its Laplacian
and the strong maximum principle force $g_{jk}=0$ for all $j,k$.
The real-analytic identity with the resulting affine modulus extends
through the positive component. Its boundary lies on the affine
complex zero hyperplane, whose complement is connected. The
component is open and relatively closed in that complement, so
continuity gives $w=|a\cdot z+b|$ on all of $\mathbb C^3$, with
$|a|=1$.

A fixed unitary affine change makes $w=|z^1|$. On a product cylinder,
the full complex Hessian controls its $(1,\bar1)$ entry. Fubini
therefore selects, for each $i$, a transverse slice whose weighted
planar energy tends to zero. Those slices converge uniformly to
$|z^1|$, since convergence is uniform on the whole cylinder.
Lemma~\ref{lem:planar-cap} gives the required contradiction.
\end{proof}

\begin{lemma}\label{lem:history}
There are $\delta,C>0$ and $C_j<\infty$ with the following property.
Suppose a solution of \eqref{eq:flow} satisfies
$$
\|R_f(t)\|_\infty+\|df(t)\|_\infty^2\le q
\quad\hbox{on }[a,b],\qquad q>0.
$$
On $[a+\delta q^{-2},b]$, the two spatial norm functions
$\|R_f(t)\|_\infty$ and $\|df(t)\|_\infty^2$ are
$Cq^3$-Lipschitz, and
\begin{equation}\label{eq:history-interpolation}
\frac{\|R_f(t)\|_\infty}{q}
\le C\left(\frac{\|df(t)\|_\infty^2}{q}\right)^{1/4},
\qquad
\|\nabla^j\Ric^{\rm Ch}(t)\|_\infty\le C_jq^{1+j/2}.
\end{equation}
\end{lemma}

\begin{proof}
Work on the periodic affine cover. At an arbitrary $(t_*,x_*)$,
use the affine scale
$\ell=(qe^{f(t_*,x_*)})^{-1/2}$ and set
$$
v(\tau,y)=f(t_*+q^{-2}\tau,x_*+\ell y)-f(t_*,x_*).
$$
Then $v(0,0)=0$, $|R_v|\le1$,
$|\nabla^0(e^{-v/2})|_E^2\le1/8$, and
$$
v_\tau=\mathcal N(v)+q^{-2}E(t_*+q^{-2}\tau),\qquad
0\le q^{-2}E\le1.
$$
The last inequality uses the original unit volume and $E\le\|R_f\|_\infty^2$.
While $|v(\tau,0)|<1$, the inverse-factor gradient bound supplies a
fixed ball on which $|v|$, $|\nabla^0v|_E$ and $|\Delta v|$ are bounded.
The exact divergence identity
\begin{align*}
\mathcal N(v)={}&-18\Delta(e^{-2v}\Delta v)
-108\diver(e^{-2v}\nabla^0v\Delta v)\\
&+45e^{-2v}(\Delta v)^2-162e^{-2v}|\nabla^0v|_E^2\Delta v
\end{align*}
therefore gives, for a radial averaging kernel $\eta_\varepsilon$,
$$
\left|\frac d{d\tau}\int\eta_\varepsilon v\right|
\le C\varepsilon^{-2},\qquad
\left|v(\tau,0)-\int\eta_\varepsilon v(\tau)\right|
\le C\varepsilon^2.
$$
The second estimate follows from the spherical mean formula and the
Laplacian bound, not a Hessian bound. Integrating from $\tau$ to zero
and choosing $\varepsilon=|\tau|^{1/4}$ gives
$|v(\tau,0)|\le C\sqrt{|\tau|}$. A first-exit argument now retains
$|v(\tau,0)|<1$ for a fixed backward time. Lemma~\ref{lem:smoothing}
on smaller cylinders gives all spatial derivatives of $v$ and
$v_\tau$, with uniform constants.

The identities
$$
(R_f)_t=-f_tR_f-3L_ff_t,\qquad
\partial_t|df|_{\gamma_f}^2
=-f_t|df|_{\gamma_f}^2+2\langle df,df_t\rangle_{\gamma_f}
$$
then bound both norm variations by $Cq^3$. Spatial suprema preserve
the resulting Lipschitz estimate. At the upper time, on a fixed spatial ball,
$\sup|\nabla^0v|_E\le C(\|df(t_*)\|_\infty^2/q)^{1/2}$ and
$\sup|\nabla^{0,3}v|_E\le C$.
The centered difference estimate for $\nabla^0v$ gives
$|\nabla^{0,2}v(0,0)|_E\le C(\|df(t_*)\|_\infty^2/q)^{1/4}$;
when the gradient norm is zero the conclusion is immediate.
This proves the first part of \eqref{eq:history-interpolation}.
The controlled affine derivatives and the two-sided bound on $e^v$
control the Chern connection and its ordered derivatives.
Undoing the metric scale gives its remaining assertion.
\end{proof}

\begin{proof}[Proof of Theorem~\ref{thm:ricci-criterion}]
Lemma~\ref{lem:ricci-gradient} and \eqref{eq:ricci-budget} imply
$$
\sup_{0<s<s_0}\int_{T-2s}^{T-s}\|df(t)\|_\infty^4\,dt<\infty.
$$
Suppose that
$\sqrt{T-t}(\|R_f(t)\|_\infty+\|df(t)\|_\infty^2)$ is unbounded.
Choose $T_i\uparrow T$ and maximize
$$
t(T_i-t)\bigl(\|R_f(t)\|_\infty+\|df(t)\|_\infty^2\bigr)^2
$$
at $t_i$, with maxima tending to infinity. Set
$K_i=\|R_f(t_i)\|_\infty+\|df(t_i)\|_\infty^2$.
Then $t_i\to T$, $K_i\to\infty$,
$t_iK_i^2\to\infty$ and $(T_i-t_i)K_i^2\to\infty$.
On their expanding time domains define
\begin{align*}
a_i(\tau)&=K_i^{-1}\bigl(
\|R_f(t_i+K_i^{-2}\tau)\|_\infty
+\|df(t_i+K_i^{-2}\tau)\|_\infty^2\bigr),\\
b_i(\tau)&=K_i^{-1}\|df(t_i+K_i^{-2}\tau)\|_\infty^2.
\end{align*}
The maximizing property gives, on each compact time interval,
$$
a_i(\tau)^2\le
\left[\left(1+\frac{\tau}{t_iK_i^2}\right)
\left(1-\frac{\tau}{(T_i-t_i)K_i^2}\right)\right]^{-1}.
$$
Lemma~\ref{lem:history} gives uniform local Lipschitz bounds.
A diagonal subsequence converges locally uniformly to globally
Lipschitz functions $a,b$ on $\mathbb R$, with
$0\le b\le a\le1$ and $a(0)=1$.
The terminal dyadic interval with $s=2(T-t_i)/3$ is centered at $t_i$;
its rescaled half-length tends to infinity. The gradient budget
therefore gives $\int_{\mathbb R}b^2<\infty$.
Global Lipschitz continuity forces $b(\tau)\to0$ as
$\tau\to-\infty$, and \eqref{eq:history-interpolation} then gives
$a(\tau)\to0$ as well.

If $a\le q$ on a finite interval, local uniform convergence gives
$a_i\le q+\varepsilon$ there. Apply Lemma~\ref{lem:history} at scale
$(q+\varepsilon)K_i$, pass to the limit and then let
$\varepsilon\downarrow0$. After delay $\delta q^{-2}$,
$a,b$ have Lipschitz constants $Cq^3$ and
$$
\frac{a-b}{q}\le C(b/q)^{1/4}.
$$
If the limiting cap holds on an entire past half-line, choose a
sufficiently early finite left endpoint for each fixed comparison.
Thus the same estimates hold on that half-line; no all-past cap for
the approximating solutions is assumed.

Let $q_j=2^{-j}$ and let $\tau_j$ be the first time $a=q_j$.
The interpolation inequality forces $b(\tau_j)\ge c q_j$ for a
fixed $c>0$. Choose $\eta>0$ small in terms of the Lipschitz
constant and $c$. On
$[\tau_j-\eta q_j^{-2},\tau_j]$, one then has
$a\ge3q_j/4$ and $b\ge c q_j/2$.
These intervals lie strictly after the first crossing of $q_j/2$
and have pairwise disjoint interiors. Each contributes at least
$\eta c^2/4$ to $\int b^2$, a contradiction.
We have proved a Type-I bound for
$\|R_f\|_\infty+\|df\|_\infty^2$.

For $\tau=T-t$, apply Lemma~\ref{lem:history} on $[t-\tau/2,t]$
with $q=C_*\tau^{-1/2}$, choosing $C_*$ large enough that
$\delta q^{-2}\le\tau/2$. For every fixed $j\ge0$ this gives
$$
\|\nabla^j\Ric^{\rm Ch}(t)\|_\infty
\le C_j(T-t)^{-1/2-j/4}.
$$
The case $j=0$ supplies the stated Ricci bound. Earlier compact
times are smooth. Compatibility with the actual singular solutions
is proved in Proposition~\ref{prop:compatibility}.
\end{proof}

On this flat conformal class the identity
$$
\mathcal R^{\rm Ch}_{i\bar jk\bar\ell}
=\tfrac13\Ric^{\rm Ch}_{i\bar j}g_{k\bar\ell},\qquad
|\mathcal R^{\rm Ch}|=|\Ric^{\rm Ch}|/\sqrt3
$$
also holds after ordered Chern differentiation, so the same rates
control the full Chern curvature tensor.

\subsection{A logarithmically improved scalar bound}

We first establish the static estimate needed to apply
Lemma~\ref{lem:history} using scalar curvature alone.

\begin{lemma}\label{lem:anchor-average}
Let $R\ge2$, $1\le q<3/2$, and let $w\ge0$ be smooth on a
neighborhood of $\overline{B_{2R}}\subset\mathbb R^6$.
If $\Delta w\le\eta$ for a continuous nonnegative function $\eta$,
then
$$
\|w\|_{L^q(B_R)}
\le C_{R,q}\left(\int_{B_{1/2}}w\,dX
+\|\eta\|_{L^1(B_{2R})}\right).
$$
\end{lemma}

\begin{proof}
Let $v\ge0$ solve $-\Delta v=\eta$ in $B_{2R}$ with zero
Dirichlet boundary values. The positive Dirichlet Green kernel is
bounded above by $C|X-Y|^{-4}$. Its integral in $X$ over $B_{1/2}$
is uniformly bounded for $Y\in B_{2R}$, so
$\int_{B_{1/2}}v\le C_R\|\eta\|_1$.
The function $z=w+v$ is nonnegative and superharmonic.

Use Green representation on $D=B_{7R/4}$ to write
$$
z=H+\int_D G_D(\,\cdot\,,Y)\,d\nu(Y),\qquad
\nu=-\Delta z\ge0,
$$
where $H\ge0$ is harmonic. Interior harmonic Harnack bounds
$\sup_{B_R}H$ by $C_R\int_{B_{1/2}}H$.
For $Y\in\overline{B_{3R/2}}$, the quantity
$\int_{B_{1/2}}G_D(X,Y)\,dX$ has a positive lower bound depending
only on $R$. Indeed it is positive and continuous in $Y$ on this
compact subset of $D$; continuity across the singularity follows
from the local integrability of $|X-Y|^{-4}$.
It follows that $\nu(B_{3R/2})\le C_R\int_{B_{1/2}}z$.
The same kernel bound gives a uniform $L^q(B_R)$ norm for
$G_D(\cdot,Y)$ when $q<3/2$. Minkowski's inequality therefore
controls the potential of this inner part of $\nu$.
For $Y\in D\setminus B_{3R/2}$ the kernel is positive harmonic
in $X\in B_{3R/2}$, so its supremum on $B_R$ is bounded by
$C_R\int_{B_{1/2}}G_D(X,Y)\,dX$. Integrating this inequality
controls the outer potential. Thus
$$
\|z\|_{L^q(B_R)}\le C_{R,q}\int_{B_{1/2}}z.
$$
Combining this with $w\le z$ and the bound for the anchor integral
of $v$ proves the result. Green representation may equivalently
be applied after smoothing the nonnegative source and passing to
the limit; all displayed bounds depend only on its $L^1$ norm.
\end{proof}

\begin{lemma}\label{lem:scalar-gradient}
For every $K<\infty$ there is $C=C(K,M,\omega_0)$ such that
every smooth real conformal factor $f$ satisfying
$\|R_f\|_{L^{3,\infty}(M,d\mu_f)}^3\le K$ obeys
$$
\|df\|_{L^\infty(\gamma_f)}^2
\le C\|R_f\|_{L^\infty(M)}.
$$
\end{lemma}

\begin{proof}
If $A=\|R_f\|_\infty>0$, partition its nonzero curvature levels
into $Ae^{-(j+1)}<|R_f|\le Ae^{-j}$, $j=0,1,\ldots$.
The weak $L^3$ bound gives
$$
\int_{\{Ae^{-(j+1)}<|R_f|\le Ae^{-j}\}}|R_f|^3\,d\mu_f
\le e^3K.
$$
On the same set, $\log(e+A/|R_f|)\ge c(j+1)$.
Summing the convergent series $\sum_{j\ge0}(j+1)^{-2}$ yields
\begin{equation}\label{eq:weak-log-budget}
\int_M\frac{|R_f|^3}{[\log(e+A/|R_f|)]^2}\,d\mu_f\le CK,
\end{equation}
where the integrand is zero at $R_f=0$.
When $A=0$, the scalar equation makes $f$ harmonic on the fixed
torus, hence constant, and the conclusion follows.

Suppose the asserted estimate fails. There are smooth factors
$f_i$ with the stated weak $L^3$ bound for which the two curvature
and gradient scales
$$
A_i=\|R_{f_i}\|_\infty>0,\qquad
B_i=\|df_i\|_{L^\infty(\gamma_{f_i})}^2>0
$$
satisfy $A_i/B_i\to0$.
Set $U_i=e^{-f_i/2}$ and choose $x_i$ where
$8|\nabla^0U_i|_E^2=B_i$. On the complete periodic affine cover put
$$
\rho_i=\frac{U_i(x_i)}{\sqrt{B_i}},\qquad
F_i(Y)=f_i(x_i+\rho_iY)-f_i(x_i),\qquad u_i=e^{-F_i/2}.
$$
The affine convention of Section~\ref{sec:equation} gives exactly
\begin{equation}\label{eq:scalar-gradient-rescaling}
u_i(0)=1,\quad |\nabla^0u_i|_E\le1/\sqrt8,\quad
|\nabla^0F_i(0)|_E=1/\sqrt2,\quad
\|R_{F_i}\|_\infty=A_i/B_i\to0.
\end{equation}

We must first bound $\rho_i$ above in order to transfer
\eqref{eq:weak-log-budget} to every fixed ball of this cover.
Without spatial dilation, the function
$U_i(x_i+X)/U_i(x_i)$ on the fixed torus has Lipschitz constant
$1/(\sqrt8\rho_i)$ and equals one at zero.
If $\rho_i$ exceeds a sufficiently large fixed constant, this
function lies between $1/2$ and $3/2$ everywhere. Its logarithm
and the scalar equation then give
$$
\|\Delta_X(f_i(x_i+\cdot)-f_i(x_i))\|_\infty
\le C\frac{A_i}{B_i\rho_i^2}.
$$
The derivative of the fixed torus Green kernel is integrable,
with local singularity $O(|X|^{-5})$. The resulting estimate
$\|\nabla^0v\|_\infty\le C\|\Delta_Xv\|_\infty$, applied to
this periodic function, gives
$1/(\sqrt2\rho_i)\le CA_i/(B_i\rho_i^2)$.
This contradicts such a fixed lower threshold for $\rho_i$ as
$A_i/B_i\to0$. Thus $\rho_i$ is bounded above.

Each fixed ball $B_R$ now meets a uniformly bounded number of
rescaled fundamental cells. On a complete cell the integral in
\eqref{eq:weak-log-budget} is unchanged: $R_{F_i}=B_i^{-1}R_{f_i}$,
the rescaled physical measure acquires the factor $B_i^3$, and
the ratio inside the logarithm is unchanged. Consequently, for
every fixed $R$,
$$
\int_{B_R}
\frac{|R_{F_i}|^3e^{3F_i}}
{[\log(e+(A_i/B_i)/|R_{F_i}|)]^2}\,\frac{dY}{8}\le C_R.
$$
The inverse-factor Lipschitz bound in
\eqref{eq:scalar-gradient-rescaling} implies
$F_i(Y)\ge-\log2-\log(1+|Y|^2)$. Define on the complete cover
$$
V_i(Y)=F_i(Y)+1+\log2+\log(1+|Y|^2)\ge1.
$$
These functions are uniformly bounded on $B_{1/2}$.
Using $|R_{F_i}|^3e^{3F_i}=27|\Delta_YF_i|^3$, the preceding
integral becomes a bound with denominator
$[\log(e+(A_i/B_i)e^{F_i}/(3|\Delta_YF_i|))]^2$.
For large $i$ we have $A_i/B_i\le1$. On $|\Delta_YF_i|\ge1$,
this logarithm is at most $C+F_i^+\le CV_i$; on the complement,
$|\Delta_YF_i|^3/V_i^2\le1$. Hence
\begin{equation}\label{eq:sublinear-scalar}
\int_{B_R}\frac{|\Delta_YF_i|^3}{V_i^2}\,dY\le C_R,
\qquad
\Delta_YV_i=g_iV_i^{2/3},\quad
\|g_i\|_{L^3(B_R)}\le C_R.
\end{equation}
Here $g_i=\Delta_YV_i/V_i^{2/3}$; the added logarithmic function
has bounded Laplacian on each fixed ball.

We next obtain every finite local moment of $V_i$.
For $0<b\le1$,
$$
\Delta_Y(V_i^b)\le b|g_i|V_i^{b-1/3}.
$$
Take $b=1/3$ and apply Lemma~\ref{lem:anchor-average} with
$q=5/4$. This gives a uniform local $L^{5/12}$ bound for $V_i$.
If an exponent $0<s\le1$ is already available on every fixed
ball, choose $b=1/3+2s/3\le1$. H\"older's inequality gives
$$
\int_{B_{2R}}|g_i|V_i^{2s/3}\,dY
\le\|g_i\|_{L^3(B_{2R})}
\left(\int_{B_{2R}}V_i^s\,dY\right)^{2/3}.
$$
The same anchor lemma improves the exponent to
$5/12+5s/6$. In particular the successive exponents
$5/12$, $55/72$, and $455/432>1$ are obtained by finitely many
applications. Bounds on $B_{2R}$ yield bounds on $B_R$;
the assertion at every fixed radius makes each step legitimate.

Once $s>1$, set $r=3s/(s+2)$, so $1<r<3$.
Equation~\eqref{eq:sublinear-scalar} bounds $\Delta_YV_i$ in
$L^r$ on each fixed ball. Let $p_i$ solve
$\Delta_Yp_i=\Delta_YV_i$ on $B_{2R}$ with zero Dirichlet data.
The constant-coefficient Poisson estimate and Sobolev embedding
give
$$
\|p_i\|_{L^{3s/2}(B_{2R})}
\le C_{R,r}\|\Delta_YV_i\|_{L^r(B_{2R})}\le C_{R,s},
\qquad \frac{2}{3s}=\frac1r-\frac13.
$$
The harmonic remainder $V_i-p_i$ is nonnegative because its
boundary values are $V_i\ge1$. Interior harmonic Harnack bounds
it on $B_R$ by its integral on $B_{1/2}$, which is bounded using
the anchor values of $V_i$ and the displayed estimate for $p_i$.
Thus $s$ improves to $3s/2$. Finitely many iterations give every
prescribed finite local moment of $V_i$, uniformly in $i$.

The functions $u_i$ have a locally uniform Lipschitz limit $u\ge0$
after passage to a subsequence. On its zero set, $V_i\to\infty$;
the moment bounds and Fatou's lemma show that this set has measure
zero. Off this null set, $F_i\to F=-2\log u$, and the higher
moments give local $L^1$ convergence of $F_i$ to $F$.
For $L>1$, the scalar equation gives
$|\Delta_YF_i|\le(A_i/(3B_i))e^L$ on $\{F_i\le L\}$.
On the complementary set, \eqref{eq:sublinear-scalar} and
H\"older's inequality yield
$$
\int_{B_R\cap\{F_i>L\}}|\Delta_YF_i|\,dY
\le C_R\left(\int_{B_R\cap\{F_i>L\}}V_i\,dY\right)^{2/3}
\le C_RL^{-2/3},
$$
where the last inequality uses the local second moment.
First let $i\to\infty$ and then $L\to\infty$.
We obtain $\Delta_YF=0$ distributionally on $\mathbb R^6$.

The smooth harmonic representative of $F$ retains the logarithmic
lower bound above. For a fixed $Y$ and a large ball of radius $R$
centered there, adding a constant of size $O_Y(\log R)$ makes
$F$ nonnegative on the ball. Differentiating the Poisson formula
and using the mean-value identity bounds
$|\nabla^0F(Y)|_E$ by $C_Y(1+\log R)/R$.
Letting $R\to\infty$ shows that $F$ is constant.
On $B_{1/2}$ the functions $u_i$ are uniformly separated from zero,
so $F=-2\log u$ there pointwise and $F(0)=0$.
It follows that $F=0$, and continuity gives $u\equiv1$.
Thus $F_i\to0$ uniformly on the anchor ball; the scalar equation
also gives $\Delta_YF_i\to0$ uniformly there.
The interior Poisson gradient estimate now implies
$\nabla^0F_i(0)\to0$, contradicting
\eqref{eq:scalar-gradient-rescaling}.
\end{proof}

\begin{proof}[Proof of Theorem~\ref{thm:scalar-criterion} and
Remark~\ref{rem:scalar-gauge}]
For each fixed time, the map
$r\mapsto r[\log(e+\sqrt{T-t}\,r)]^{1/3}$ is a continuous
strictly increasing bijection of $[0,\infty)$.
The definition of the weak $L^3$ quasi-norm therefore identifies
the cube of the norm in \eqref{eq:scalar-log} exactly with
$$
\sup_{\lambda>0}\lambda^3\log(e+\sqrt{T-t}\lambda)
\mu_{f(t)}\{|R_f(t)|>\lambda\}.
$$
It suffices to prove the conclusion under
\eqref{eq:scalar-gauge}. Let $K$ be a finite upper bound there.
Since $h\ge1$, Lemma~\ref{lem:scalar-gradient} gives, for every
$t_0<t<T$,
$$
\|R_f(t)\|_{L^{3,\infty}(d\mu_f(t))}^3\le K,\qquad
\|df(t)\|_\infty^2\le C_K\|R_f(t)\|_\infty.
$$

If the solution were Chern Type-II, choose $\sigma_i\uparrow T$
with $Q_i=\|R_f(\sigma_i)\|_\infty\to\infty$ and
$\sqrt{T-\sigma_i}Q_i\to\infty$.
Let $b_i$ be the first time at which $\|R_f\|_\infty=Q_i$.
Continuity and smoothness on compact preterminal intervals give
$$
b_i\le\sigma_i,\qquad b_i\to T,\qquad
\|R_f(t)\|_\infty\le Q_i\quad(0\le t\le b_i).
$$
On the part of this past lying after $t_0$, the gradient estimate
gives $\|R_f\|_\infty+\|df\|_\infty^2\le(1+C_K)Q_i$.
For large $i$ this cap holds throughout a past interval of twice
the delay required by Lemma~\ref{lem:history} at scale
$(1+C_K)Q_i$. The case $j=1$ of that lemma yields
$$
\|dR_f(b_i)\|_{L^\infty(\gamma_f(b_i))}
\le C_KQ_i^{3/2}.
$$
Here $R_f$ is the trace of first Chern--Ricci curvature, and the
Chern connection preserves the metric, so its differentiated
trace is controlled by $\nabla\Ric^{\rm Ch}$.

Choose $x_i$ with $|R_f(b_i,x_i)|=Q_i$, and set at time $b_i$
$$
U_i=e^{-f(b_i)/2},\qquad r_i=U_i(x_i)/\sqrt{Q_i}.
$$
The gradient estimate gives
$|\nabla^0U_i|_E\le\sqrt{C_KQ_i/8}$.
Consequently, on a fixed sufficiently small multiple of
$B_0(x_i,r_i)$ in the affine cover,
$U_i/U_i(x_i)$ lies between $1/2$ and $3/2$.
Unit volume and the unit background volume give a point where
$U_i\le1$, so
$$
r_i\le Q_i^{-1/2}+\sqrt{C_K/4}\,\operatorname{diam}_{\gamma_0}(M).
$$
This uniformly bounds $r_i$ for large $i$.
We may therefore fix $\eta>0$, independent of $i$, such that
$B_0(x_i,\eta r_i)$ embeds in the torus and has the preceding
two-sided bound. A straight affine segment in this ball has
$\gamma_f(b_i)$-length at most $C\eta Q_i^{-1/2}$.
The bound for $dR_f$ then gives
$|R_f(b_i)|>Q_i/2$ throughout the ball after reducing $\eta$.
Its physical volume satisfies
$$
\mu_{f(b_i)}\bigl(B_0(x_i,\eta r_i)\bigr)
\ge \frac{|B_1|}{8}(3/2)^{-6}\eta^6
U_i(x_i)^{-6}r_i^6=c_KQ_i^{-3},
$$
with $c_K>0$.

Testing \eqref{eq:scalar-gauge} at time $b_i$ and level
$\lambda=Q_i/2$ now gives
$$
K\ge\frac{c_K}{8}
h\!\left(\tfrac12\sqrt{T-b_i}\,Q_i\right).
$$
But $\sqrt{T-b_i}Q_i\ge\sqrt{T-\sigma_i}Q_i\to\infty$,
contradicting $h(s)\to\infty$. This proves the asserted terminal
Type-I bound. Smoothness controls all earlier compact time
intervals and completes both statements.
\end{proof}

\section{Radial shrinking solutions}\label{sec:shrinking}

We now work in complex dimension three. The unnormalized operator
is \eqref{eq:N}, with the Euclidean Laplacian in the affine
coordinates fixed in Section~\ref{sec:equation}. The profiles below
are smooth at the origin and approach a logarithmic cone at
infinity. Put $\kappa=\sqrt5-1$.

\begin{lemma}\label{lem:profile}
There are sequences $h_j\to\infty$ and $k_j\to\kappa$, with $k_j>0$ and
$|k_j-\kappa|=O(h_j^{-1/2})$, and smooth radial functions $F_j$ on $\mathbb R^6$ such that
$$
F_j(0)=0,\qquad \Delta F_j(0)=h_j,\qquad F_j'(r)>0\quad(r>0),\qquad\Delta F_j(r)\ge c_j(1+r^2)^{-1}>0.
$$
If
$$
b_j=\frac1{2(k_j+2)},\qquad a_j=k_jb_j,
$$
then
$$
-18e^{-3F_j}\Delta(e^{F_j}\Delta F_j)
-9e^{-2F_j}(\Delta F_j)^2=-a_j+b_jrF_j'(r).
$$
There are real constants $c_{\infty,j}$ and $c_{j,\ell}$ such that, for every pair of nonnegative integers $N,d$,
$$
\frac{d^d}{dr^d}\left[
F_j(r)-k_j\log r-c_{\infty,j}
-\sum_{\ell=1}^N c_{j,\ell}r^{-2(k_j+2)\ell}
\right]
=O\!\left(r^{-2(k_j+2)(N+1)-d}\right).
$$
Consequently, for every $T>0$, the function
$$
f(t,x)=a_j\log(T-t)+F_j\!\left(\frac{|x|}{(T-t)^{b_j}}\right)
$$
solves the displayed evolution equation on $\mathbb R^6$ for $t<T$.
\end{lemma}

The proof is given in Appendix~\ref{sec:profile}. It includes the
two linear connection calculations needed to match the regular
origin to the end at infinity. In particular, the profile selection
does not rely on numerical shooting.

Fix one of these profiles and omit its index. Write
$$
\gamma_F=\tfrac12e^{F(|y|)}\sum_{j=1}^6(dy^j)^2.
$$
The factor $\tfrac12$ agrees with the affine real metric convention
in Section~\ref{sec:equation}. Since $F$ is smooth
at zero and $F(r)=k\log r+c_\infty+O(r^{-2(k+2)})$ with all
finite derivatives, $\gamma_F$ is complete: its radial length to
infinity is comparable to $\int_1^\infty r^{k/2}\,dr$.
Its Chern scalar curvature is
$R_F=-3e^{-F}\Delta F$. Its Chern curvature is bounded and decays
as $O(r^{-k-2})$. At the origin $R_F(0)=-3h<0$, so
$$
0<\sup_{\mathbb R^6}|R_F|<\infty.
$$
For $b=1/[2(k+2)]$ and $a=kb$, the identity $a+2b=1/2$
converts the solution in Lemma~\ref{lem:profile} into the complete
shrinking family
$$
\gamma_\infty(q)=\tfrac12(-q)^a
\exp\!\left(F\!\left(\frac{|y|}{(-q)^b}\right)\right)
\sum_{j=1}^6(dy^j)^2,\qquad q<0.
$$
It solves $f_t=\mathcal N(f)$. The next section realizes its
asymptotics on the fixed compact torus.

\section{Localization on the torus}\label{sec:localization}

\begin{proposition}\label{prop:localized}
For every flat complex three-torus of volume one and every $p\in M$, there are a
profile $F$ from Lemma~\ref{lem:profile}, positive numbers
$T,r_*,\delta,R_0$, and a smooth solution
$\widetilde f_s=\mathcal N(\widetilde f)$ on $[0,T)\times M$
with the following properties. Write $\sigma=T-s$,
$b=1/[2(k+2)]$ and $a=kb$.
In a fixed affine chart centered at $p$,
$$
\widetilde f(s,X)=a\log\sigma+F(|X|/\sigma^b)+v(s,X)
\quad\hbox{on }B_{r_*}.
$$
For every fixed $R<\infty$ and $0\le j\le10$,
\begin{equation}\label{eq:core-output}
\sup_{|y|\le R}|\nabla^{0,j}_{y}v(s,\sigma^by)|
\le C_R\sigma^\delta
\end{equation}
for all sufficiently late $s<T$. For $0\le j\le10$,
\begin{equation}\label{eq:annular-output}
|\nabla^{0,j}_{X}v(s,X)|\le C|X|^{\delta/b-j}
\quad\hbox{if }R_0\sigma^b\le|X|\le r_*.
\end{equation}
The $C^{10}$ norm of $\widetilde f$ on
$M\setminus B_{r_*/2}$ is uniformly bounded. Its volume satisfies
$$
0<\inf_{s<T}\int_Me^{3\widetilde f(s)}d\mu_0
\le\sup_{s<T}\int_Me^{3\widetilde f(s)}d\mu_0<\infty.
$$
All constants are independent of the preterminal time.
\end{proposition}

\begin{proof}
We first construct a solution $\widetilde f_s=\mathcal N(\widetilde f)$ and then impose the unit-volume normalization.
Lemma \ref{lem:profile} supplies a smooth radial $F$ on $\mathbb R^6$ and a number $k>0$ such that, with $b=1/(2(k+2))$ and $a=bk$,
$$
\mathcal N(F)=-a+brF_r,\quad F(0)=0,\quad
F(r)=k\log r+c_\infty+O(r^{-2k-4}).
$$
The expansion can be differentiated to every prescribed finite order. For a small $T>0$, the exact local solution is
$$
f_*(s,X)=a\log(T-s)+F(|X|/(T-s)^b),\qquad 0\le s<T.
$$

Choose $r_*>0$ so that $\overline B_{10r_*}$ lies in the affine chart centered at $p$. Choose $f_\infty$ smooth on $M\setminus\{p\}$ and equal to $k\log|X|+c_\infty$ on $B_{8r_*}\setminus\{0\}$. Let $\chi_{\rm app}$ equal one on $B_{6r_*}$ and vanish outside $B_{7r_*}$. Define
$$
f_{\rm app}=f_\infty+\chi_{\rm app}(f_*-k\log|X|-c_\infty).
$$
It is defined by $f_*$ at the origin. On each fixed exterior region its spatial derivatives and its residual $\partial_sf_{\rm app}-\mathcal N(f_{\rm app})$ are uniformly bounded as $T\downarrow0$ and $s\uparrow T$. Indeed the profile has an expansion in powers $r^{-2(k+2)j}$, which becomes an expansion in integer powers of $T-s$ on a fixed physical annulus. Its differentiated remainders justify these bounds, including the time derivative.

Put $\sigma=T-s$, $\tau=-\log\sigma$, $y=X/\sigma^b$, and $v=\widetilde f-f_{\rm app}$. Choose $\chi$ equal to one on $B_{2r_*}$ and supported in $B_{3r_*}$. The function
$$
w(\tau,y)=\chi(\sigma^by)v(s,\sigma^by)
$$
is extended by zero to $\mathbb R^6$. On this support $f_{\rm app}=f_*$. Consequently
\begin{equation}\label{eq:w}
 w_\tau=\mathcal Aw+Q(w)+\mathcal E,\qquad
 \mathcal A=\mathcal N'_F-by\cdot\nabla^0,
\end{equation}
where $\mathcal N'_F$ is the linearization of $\mathcal N$ at $F$,
$Q(w)=\mathcal N(F+w)-\mathcal N(F)-\mathcal N'_Fw$, and
\begin{equation}\label{eq:cutoff}
\mathcal E=\chi Q(v)-Q(\chi v)-[\mathcal N'_F,\chi]v.
\end{equation}
The expression $v$ in this formula is in rescaled coordinates. The terms $\chi_\tau=-by\cdot\nabla^0\chi$ and the transport commutator cancel exactly. Thus \eqref{eq:cutoff} is supported in $2r_*\le|X|\le3r_*$.

\medskip\noindent\textbf{1. Estimates in the rescaled ball.}
Lemma \ref{lem:linear} realizes $\mathcal A$ in
$$
L^2(e^{-2\Phi}d\mu_{g_F}),\qquad
g_F=e^F\sum_{j=1}^6(dy^j)^2,\quad \Phi=\alpha\rho^{4/3}
$$
outside a compact set, where $\rho$ is proportional to $r^{(k+2)/2}$ and $\alpha>0$ is sufficiently small. Its resolvent is compact. Choose $0<\delta<1/16$ such that the line $\Re\lambda=-\delta$ avoids its spectrum. Let $\Pi_+$ project onto the finite-dimensional generalized eigenspaces to its right, and let $\Pi_-=I-\Pi_+$. Here $\nabla^F$ denotes the Chern connection of $g_F$, and $\nabla^{F,j}=(\nabla^F)^j$. The spaces $\mathcal G_j$ of Lemma \ref{lem:linear} have norms
$$
\|w\|_{\mathcal G_j}^2=
\sum_{\ell=0}^j\|\langle\rho\rangle^{(j-\ell)/3}
\nabla^{F,\ell}w\|_{L^2(e^{-2\Phi}d\mu_{g_F})}^2.
$$
The stable equation has maximal regularity with base, trace and domain $\mathcal G_{12},\mathcal G_{14},\mathcal G_{16}$, respectively, on every backward unit $\tau$-interval, with exponential weight $e^{\delta\tau}$ and a constant independent of the terminal time. We use exactly this causal uniform-local norm, denoted $\mathfrak X_\delta$ in that lemma.

If $\sum_{j\le8}\|\nabla^{F,j}w\|_\infty$ is sufficiently small, Lemma \ref{lem:tame} gives
\begin{equation}\label{eq:tame-use}
\begin{split}
\|Q(w)\|_{\mathcal G_{12}}
&\le C\left(\sum_{j\le8}\|\nabla^{F,j}w\|_\infty\right)
\|w\|_{\mathcal G_{16}},\\
\|\Pi_+Q(w)\|&\le
C\left(\sum_{j\le8}\|\nabla^{F,j}w\|_\infty\right)
\|w\|_{\mathcal G_{14}}.
\end{split}
\end{equation}
The second estimate is instantaneous. In the first, after $\ell\le12$ derivatives, the total derivative order is at most sixteen. All factors other than one of highest order have order at most eight; the weights satisfy $12-\ell\le16-j$ for the chosen derivative order $j$. The same argument at order zero gives $\|Q(w)\|_{\mathcal G_0}\le C\|w\|_{C^8(g_F)}\|w\|_{\mathcal G_4}$, proving the second inequality. An algebra property of the weakly weighted space is not being assumed.

\medskip\noindent\textbf{2. Estimates on physical annuli and the exterior.}
Choose a large fixed $R$. For each dyadic radius $r\le r_*$, use the central annulus $r/2<|X|<2r$ and the enlarged annulus $r/4<|X|<4r$, starting at
$$
s_r=\max\{0,T-(4r/R)^{2(k+2)}\}.
$$
Set $u(\theta,z)=v(s,rz)$, where $\theta=(s-s_r)/r^{2(k+2)}$.
The active time interval has length at most $4^{2(k+2)}R^{-2(k+2)}$.
We propagate the $L^\infty H^{14}\cap L^2H^{16}$ bound on the central annulus, $\|u\|_{L^\infty H^{14}\cap L^2H^{16}}\le\Lambda\eta r^{\delta/b}$, where $\Lambda$ is fixed after $R$
and $\eta>0$ is small.

The differentiated tail gives uniformly bounded coefficients and uniform strong parabolicity in these variables. Let $A_r$ be the linearization at $f_*$ expressed in these annular coordinates, and let $Q_r$ be the remainder after subtracting the linear part. Then $u_\theta=A_ru+Q_r(u)$ and
$$
\|Q_r(u)\|_{H^{12}}\le C\|u\|_{H^{14}}\|u\|_{H^{16}}.
$$
For a cutoff $\zeta$ supported in the enlarged annulus and equal to one on the central annulus,
$$
(\zeta u)_\theta=A_r(\zeta u)+[\zeta,A_r]u+\zeta Q_r(u).
$$
The commutator has order three. Integrating by parts yields
$$
\frac d{d\theta}\|\zeta u\|_{H^{14}}^2+c\|\zeta u\|_{H^{16}}^2
\le C\|\zeta u\|_{H^{14}}^2+C\|[\zeta,A_r]u+\zeta Q_r(u)\|_{H^{12}}^2.
$$
On an interval of length $\ell\le1$, interpolation gives
$$
\|[\zeta,A_r]u\|_{L^2H^{12}}
\le C\|u\|_{L^2H^{15}}
\le C\ell^{1/4}\|u\|_{L^\infty H^{14}\cap L^2H^{16}}.
$$
Hence the central bound improves to
\begin{equation}\label{eq:annulus}
\|u\|_{L^\infty H^{14}\cap L^2H^{16}}
\le C_R\eta r^{\delta/b}+C R^{-(k+2)/2}\Lambda\eta r^{\delta/b}
   +C \Lambda^2\eta^2r^{2\delta/b}.
\end{equation}

Here is the complete overlap and initial-data justification. The enlarged annulus is covered by the central annuli of radii $r/2,r,2r$. The $2r$ interval has already started. Before the $r/2$ interval starts, its overlap lies in $R/16\le|y|\le R$, on a bounded $\tau$-interval; there $\sigma^\delta$ is comparable to $r^{\delta/b}$. Thus the core trace and local time-integrated domain norms provide $C_R\eta r^{\delta/b}$. At a positive starting time $s_r$, the entire enlarged annulus lies in $|y|\le R$, giving the same initial bound. At $s_r=0$ its initial perturbation will be zero by the support choice below. The constants in these comparisons are fixed once $R$ is fixed. This accounts for all lateral fluxes in the displayed commutator.

For the fixed exterior choose a cutoff $\zeta_0$ vanishing on $B_{r_*/4}$ and equal to one outside $B_{r_*/2}$. Work with $\zeta_0v$ on the whole torus. Its cutoff commutator is supported in a buffer controlled by annuli of fixed radii comparable to $r_*$. The approximate coefficients are uniformly bounded there; extend them smoothly and uniformly elliptically inside the region where $\zeta_0v=0$. Its equation contains this order-three commutator, the nonlinear remainder, and the bounded residual of $f_{\rm app}$. The same energy estimate, in physical time, gives
\begin{equation}\label{eq:exterior}
\|\zeta_0v\|_{L^\infty H^{14}\cap L^2H^{16}}
\le C_{r_*}\{T^{1/2}+T^{1/4}\Lambda\eta+(\Lambda\eta)^2\}
\end{equation}
under the exterior bootstrap bound $\Lambda\eta$. The initial exterior perturbation is zero. The largest enlarged annuli use this exterior norm on their overlap. All conversion factors depend only on the already fixed $r_*$, $\delta$ and $b$ and are included when choosing $R$. Thus there is no unprescribed boundary value in \eqref{eq:annulus} or \eqref{eq:exterior}.

Choose $R$ large so the linear coefficient in \eqref{eq:annulus} is small, then $\Lambda$ large compared with the core restriction constants, then $\eta$ small. Finally take $T$ small. Both outer bounds improve by a fixed factor. At every $s_1<T$ there are only finitely many active dyadic radii, all bounded below by a positive multiple of $R(T-s_1)^b$; the argument therefore loses no endpoint at radius zero.

Sobolev embedding gives $|\nabla^{0,j}_{X}v|\le C\Lambda\eta r^{\delta/b-j}$ for $j\le10$. Passing to $y$ and then to $g_F$ gives
$$
|\nabla^{F,j}v|\le C\Lambda\eta r^{\delta/b}\rho^{-j},\qquad j\le8,
$$
outside the fixed rescaled core. The physical cutoff has intrinsic derivatives $O(\rho^{-j})$ there. Together with the $\mathcal G_{14}$ core trace, this supplies the uniform smallness needed in \eqref{eq:tame-use}.

On the support of \eqref{eq:cutoff}, $\rho\asymp e^{\tau/4}$. All coordinate, derivative and time conversions cost only $e^{C\tau}$, whereas the weight is at most $e^{-c e^{\tau/3}}$. The preceding outer $H^{14}$ trace and $L^2H^{16}$ estimates therefore give
\begin{equation}\label{eq:error}
\begin{split}
\|\mathcal E\|_{L^2([\max\{-\log T,\tau-1\},\tau];\mathcal G_{12})}
&\le C\Lambda\eta e^{C\tau-c e^{\tau/3}},\\
\|\Pi_+\mathcal E(\tau)\|&\le C\Lambda\eta e^{C\tau-c e^{\tau/3}}.
\end{split}
\end{equation}
The second line only requires four spatial derivatives, hence the trace suffices. The first uses the time-integrated sixteenth derivatives. Taking the initial $\tau_0=-\log T$ large makes both weighted error norms as small as required for all later times.

\medskip\noindent\textbf{3. Selection of a genuine solution.}
Let $U=\Ran\Pi_+$ be the real invariant subspace, taking conjugate eigenvalues together. The map $\mathcal J$ and all parameters below are real. By Lemma \ref{lem:linear}, there is a linear map $\mathcal J:U\to C_c^\infty(\mathbb R^6)$ with $\Pi_+\mathcal J=I$, support in $|y|\le R_{\rm mode}$, and $\|\Pi_-\mathcal J\|_{U\to\mathcal G_{14}}\le\beta$, where $\beta>0$ is as small as prescribed. Namely, cut off the finitely many generalized eigenfunctions in their graph norms and invert the resulting finite-dimensional projected map.

All eigenvalues of $\mathcal A|_U+\delta I$ have positive real part. Give $U$ the norm associated with
$$
G=\int_0^\infty e^{-(\mathcal A|_U+\delta I)^*q}
                  e^{-(\mathcal A|_U+\delta I)q}\,dq.
$$
It satisfies $(\mathcal A|_U+\delta I)^*G+G(\mathcal A|_U+\delta I)=I$. For $|\xi|_G\le1$, choose the unnormalized initial data
\begin{equation}\label{eq:initial}
\widetilde f_0(X)=f_{\rm app}(0,X)+\eta T^\delta(\mathcal J\xi)(X/T^b),
\qquad
f_0=\widetilde f_0-\tfrac13\log\int_Me^{3\widetilde f_0}d\mu_0.
\end{equation}
The perturbation is extended by zero. These are smooth strictly positive conformal metrics, and the second expression has volume one.

Here is the order of all choices. Fix the physical chart and radii, the profile, the weight, the spectral line and the fixed geometric norm conversions. Choose $\beta$ small and then $R_{\rm mode}$. Choose $R>16R_{\rm mode}$ as in \eqref{eq:annulus}, then $\Lambda$, then $\eta$ small for the nonlinear and exit estimates. Finally choose $\tau_0$ large so $RT^b<r_*/4$, \eqref{eq:exterior} improves, and \eqref{eq:error} is sufficiently small. If an annulus starts at zero, then $r/4\ge RT^b/16>R_{\rm mode}T^b$, so its initial perturbation is exactly zero. The fixed rescaled core always lies where $\chi=1$.

Bootstrap the outer bounds above together with
$$
|e^{\delta\tau}\Pi_+w(\tau)|_G\le\eta,
\qquad \|\Pi_-w\|_{\mathfrak X_\delta}\le\eta.
$$
On the finite-dimensional component the pointwise bound controls every spatial graph norm on a unit time interval. Its time derivative follows from \eqref{eq:w}. Applying the stable estimate in Lemma \ref{lem:linear}, \eqref{eq:tame-use}, and \eqref{eq:error}, and absorbing the small coefficient, gives
$$
\|\Pi_-w\|_{\mathfrak X_\delta}
\le C\beta\eta+C(C_R+\Lambda)\eta^2+o_{\tau_0}(\eta)\le\eta/4.
$$
The outer estimates also improve strictly. The initial stable trace and outer norms satisfy their bounds strictly. Thus the only possible first exit is through $|e^{\delta\tau}\Pi_+w|_G=\eta$.

Write $z=e^{\delta\tau}\Pi_+w$. At this boundary, its equation and the instantaneous part of \eqref{eq:tame-use} imply
$$
\frac d{d\tau}|z|_G^2
=|z|^2+2\langle Gz,e^{\delta\tau}\Pi_+(Q(w)+\mathcal E)\rangle
\ge c\eta^2>0.
$$
The exit is strictly outward. Suppose every parameter in the closed unit $G$-ball exits at a finite $\tau$. Uniqueness, continuous dependence and strict transversality make the first exit point continuous in the parameter; compactness supplies a finite bound on these exit times. The map sending $\xi$ to $z/\eta$ at its first exit would retract the ball onto its sphere, with the identity on the boundary. This contradicts the no-retraction theorem~\cite[Chapter~1]{Deimling85}. If $U=\{0\}$, the stable and outer improvements alone give a nonexiting solution. Thus in either case a parameter does not exit.

For completeness the actual PDE cannot terminate at $S<T$ while these bounds hold. Since $T-s\ge T-S>0$, the rescaling and all weighted restrictions on the needed fixed balls are nonsingular. The core and outer bounds give a uniform physical $H^{14}(M)$ bound for $\widetilde f$ up to $S$. Lemma~\ref{lem:Sobolev-continuation} extends the actual solution across $S$. This rules out a premature PDE endpoint and justifies the exit construction at every finite $\tau$.

The selected solution therefore exists smoothly for all $s<T$, and at the chosen point
\begin{equation}\label{eq:center}
\widetilde f(s,p)=a\log(T-s)+O(\eta(T-s)^\delta)\longrightarrow-\infty.
\end{equation}
It cannot extend positively across $T$, so its unnormalized maximal time is exactly $T$.

\medskip\noindent\textbf{4. The original normalization and time.}
Set $V(s)=\int_Me^{3\widetilde f(s)}d\mu_0$. Here $E(\widetilde f(s))=9\int_Me^{\widetilde f}(L_0\widetilde f)^2d\mu_0$. Integrating the unnormalized equation gives
$$
V'=-3E(\widetilde f(s))\le0.
$$
The fixed exterior estimate keeps $\widetilde f$ bounded below on a fixed nonempty open set outside $B_{7r_*}$. Consequently $0<v_*\le V(s)\le V(0)$ for $s<T$. Define
\begin{equation}\label{eq:time}
t(s)=\int_0^sV(q)^{-2/3}dq,\qquad
f(t(s))=\widetilde f(s)-\tfrac13\log V(s).
\end{equation}
The constant-shift laws are $\mathcal N(f)=V(s)^{2/3}\mathcal N(\widetilde f(s))$ and $E(f)=V(s)^{-1/3}E(\widetilde f(s))$. Differentiating \eqref{eq:time} gives exactly \eqref{eq:flow}. Its terminal time obeys
$$
0<\Theta=\int_0^T V(q)^{-2/3}dq\le Tv_*^{-2/3}<\infty.
$$
Equation \eqref{eq:center} and the two volume bounds give
$e^{f(t,p)}\to0$. A positive smooth extension would have a positive
lower bound at $(\Theta,p)$, a contradiction. Uniqueness in
Proposition~\ref{prop:existence} shows that this is the maximal
normalized solution from \eqref{eq:initial}.

The weighted trace bound for $w$ restricts to an ordinary $H^{14}$
bound on every fixed rescaled ball, where the weight and its reciprocal
are bounded. Since $\chi(\sigma^by)=1$ there for late time, Sobolev
embedding gives \eqref{eq:core-output}. The dyadic annular estimate
and the exterior trace give \eqref{eq:annular-output} and the stated
exterior $C^{10}$ bound. These are bounds for the actual selected
solution, so all assertions of the proposition follow.
\end{proof}

\begin{lemma}\label{lem:alltime}
Fix the profile, the physical radius and an exterior function allowed
in Proposition~\ref{prop:localized}. Its selected solution can be
chosen with parameters $R,\Lambda,\eta,T$ in the order in that proof
so that, for every $0\le s<T$ and $0\le j\le10$,
\begin{align}
|\nabla_y^{0,j}v(s,\sigma^by)|
&\le C_R\eta\sigma^\delta
&&(|y|\le2R),\label{eq:allcore}\\
|\nabla_X^{0,j}v(s,X)|
&\le C\Lambda\eta r^{\delta/b-j}
&&(2R\sigma^b\le r\le r_*/2).\label{eq:allannulus}
\end{align}
On the fixed exterior,
\begin{equation}\label{eq:allexterior}
\sup_{s<T}\|v(s)\|_{C^2(M\setminus B_{r_*/2})}
\le C_{r_*,f_\infty}
\{T^{1/2}+T^{1/4}\Lambda\eta+(\Lambda\eta)^2\}.
\end{equation}
\end{lemma}

\begin{proof}
The selected unstable parameter stays in the fixed unit ball.
Its bound and the stable trace in $\mathcal G_{14}$ hold from
$\tau_0=-\log T$, including the initial trace. Restriction to the
fixed rescaled ball $B_{2R}$ and $H^{14}\hookrightarrow C^{10}$
give \eqref{eq:allcore}.

If a physical annulus has started, its improved estimate
in \eqref{eq:annulus} and rescaling give \eqref{eq:allannulus}.
Before it starts, its intersection with the indicated region is
covered by the core, as in the overlap argument of that proof.
Since $\sigma$ and $r^{1/b}$ are then comparable with constants
depending on $R$, the same estimate follows from
\eqref{eq:allcore}. The active central annuli cover the remaining
region with bounded overlap. At the initial time their perturbation
is either the core trace or zero by its support.
Finally \eqref{eq:exterior} and Sobolev embedding give
\eqref{eq:allexterior}. The term $(\Lambda\eta)^2$ is made small
by choosing $\eta$ before $T$.
\end{proof}

\section{Small positive and subcritical curvature throughout the flow}
\label{sec:smallness}

We choose the exterior function in Proposition~\ref{prop:localized}
so that the sign of scalar curvature can be controlled through
the whole lifetime. We write $R_f^+=\max\{R_f,0\}$ for its positive part.

\begin{proposition}\label{prop:smallness}
For a fixed sufficiently large profile of Lemma~\ref{lem:profile},
the localized solution can be chosen to satisfy
\eqref{eq:smallpositive}--\eqref{eq:smallsubcritical} and
$0<\Theta<\tau_*$ for each prescribed $p_0,\varepsilon,\tau_*$
in Corollary~\ref{cor:smallness}.
\end{proposition}

\begin{proof}
It suffices to treat $2\le p_0<3$: on the normalized probability
measure, every smaller exponent is bounded by the $L^2$ norm.
Fix a physical radius $R_{\rm ext}$ with its doubled ball in the
affine chart at $p$, and a cutoff $\zeta$ equal to one in
$B_{R_{\rm ext}}$ and zero outside $B_{2R_{\rm ext}}$.
Choose a smooth radial cutoff $\chi$ equal to one for $r\le8$
and zero for $r\ge9$, with $0\le\chi\le1$.
For $9r_*<R_{\rm ext}$ define the nonpositive radial potential
\begin{equation}\label{eq:newton}
v_{r_*}(r)=
-\int_r^\infty s^{-5}
\int_0^s4k t^3\chi(t/r_*)\,dt\,ds.
\end{equation}
Radial differentiation gives
$\Delta_Xv_{r_*}=4kr^{-2}\chi(r/r_*)\ge0$.
For $r<8r_*$ it equals $k\log r+C_{r_*}$, while for $r>9r_*$
it equals
$$
-\frac1{4r^4}\int_0^\infty4kt^3\chi(t/r_*)\,dt.
$$
The integral is a fixed positive constant times $r_*^4$.
Scaling \eqref{eq:newton} gives
$C_{r_*}=C_1-k\log r_*$ for a fixed $C_1$.
The radial flux at zero vanishes, so there is no Dirac mass there.

Use the allowed exterior
$$
f_\infty=c_\infty-C_{r_*}+\zeta v_{r_*}.
$$
It equals $k\log r+c_\infty$ on $B_{8r_*}\setminus\{p\}$.
All negative Laplacian lies in the fixed far annulus
$B_{2R_{\rm ext}}\setminus B_{R_{\rm ext}}$, and
\begin{equation}\label{eq:exterior-negative}
\|[\Delta_X(\zeta v_{r_*})]^-\|_\infty\le Cr_*^4.
\end{equation}
Indeed on this annulus the potential is a multiple of $r_*^4r^{-4}$,
and every differentiated cutoff has a fixed bound.

Let $f_{\rm app}$ be the approximation in the localization proof.
After subtracting the constant $c_\infty-C_{r_*}$ it is the exact
core $a\log\sigma+F(r/\sigma^b)-(c_\infty-C_{r_*})$,
interpolated with $v_{r_*}$ on $6r_*\le r\le7r_*$, and
equal to $\zeta v_{r_*}$ outside.
For $0\le j\le2$, the differentiated tail gives
$$
|\nabla_X^{0,j}(a\log\sigma+F(r/\sigma^b)-k\log r-c_\infty)|
\le C\sigma r^{-2(k+2)-j}\qquad(r\asymp r_*).
$$
As $F'>0$, the maximum of the shifted core on $B_{7r_*}$ is
its boundary value, which tends to
$v_{r_*}(7r_*)=v_1(7)<0$ as $T\downarrow0$.
Consequently
$$
f_{\rm app}-(c_\infty-C_{r_*})\le0
\quad(0\le s<T).
$$
On $M\setminus B_{r_*/2}$ the same shifted approximation is
bounded below by a constant depending only on $F$ and the fixed
cutoffs. Its negative Laplacian is bounded by $Cr_*^4$.
On the transition annulus, the Laplacian $4kr^{-2}$ is changed
by at most $CT r_*^{-2(k+2)-2}$ and therefore remains positive
when $T$ is chosen last.

For the full affine Hessian, uniformly in time,
\begin{equation}\label{eq:approx-subcritical}
\|\nabla^{0,2}(f_{\rm app}-(c_\infty-C_{r_*}))\|_{L^{p_0}(d\mu_0)}
\le C_{p_0}r_*^{6/p_0-2}.
\end{equation}
In the core the Hessian is at most
$C(\sigma^{2b}+r^2)^{-1}$, whose $p_0$th power integral is at most
$C_{p_0}r_*^{6-2p_0}$. The transition has the same bound.
The potential has Hessian at most $Cr^{-2}$ for $7r_*\le r\le9r_*$,
and at most $Cr_*^4r^{-6}$ for $9r_*\le r\le R_{\rm ext}$.
The latter integral is bounded by
$Cr_*^{4p_0}\int_{9r_*}^{R_{\rm ext}}r^{5-6p_0}\,dr
\le C_{p_0}r_*^{6-2p_0}$, since $p_0>1$.
The far cutoff contributes $O(r_*^4)$.

Apply Proposition~\ref{prop:localized} with this exterior.
Lemma~\ref{lem:alltime} gives the actual perturbation.
Since $\Delta F\ge c_F(1+r^2)^{-1}$, its core Laplacian is
positive at every time after choosing $T$ small in
\eqref{eq:allcore}. On active annuli the positive profile
Laplacian is at least $c_Fr^{-2}$, whereas the error is at most
$C\Lambda\eta r^{\delta/b-2}$. Choose $\eta$ small first.
It follows that $\Delta_X\widetilde f>0$ on $B_{r_*/2}$
through the whole interval.

The integrated Hessian error from \eqref{eq:allcore}--\eqref{eq:allexterior}
is bounded by
\begin{align*}
C_{p_0,R}\eta T^{\delta+b(6/p_0-2)}
+C_{p_0}\Lambda\eta r_*^{\delta/b+6/p_0-2}\\
+C_{r_*,f_\infty}
\{T^{1/2}+T^{1/4}\Lambda\eta+(\Lambda\eta)^2\}.
\end{align*}
It can be made arbitrarily small, by $\eta$ first and $T$ last.
The same choices ensure $|v|\le1$ globally and make the exterior
$C^2$ error smaller than any fixed positive target.
Thus, after the constant shift
$\bar f=\widetilde f-(c_\infty-C_{r_*})$, one has
$\bar f\le1$ everywhere, a lower bound independent of $r_*$
off $B_{r_*/2}$, and
$$
\sup_s\|R_{\bar f}^+\|_\infty\le C_F(r_*^4+o(1)).
$$
The term $o(1)$ here is made small for the already fixed $r_*$
by $\eta$ and $T$. Inside $B_{r_*/2}$ the positive scalar part
is zero.

The shift changes the unnormalized time. Precisely, set
$$
\bar s=e^{-2(c_\infty-C_{r_*})}s,\qquad
\bar T=e^{-2(c_\infty-C_{r_*})}T .
$$
Then $\bar f(\bar s)$ solves $\bar f_{\bar s}=\mathcal N(\bar f)$.
Write $\bar V(\bar s)=\int_Me^{3\bar f(\bar s)}d\mu_0$. Its initial volume tends to one as $r_*\downarrow0$ and the
remaining errors tend to zero. To see this directly, use
$v_{r_*}\le0$ and
$$
\int_M|\zeta v_{r_*}|\,d\mu_0\le Cr_*^4.
$$
The inner contribution is $O(r_*^6)$ by scaling, and the Newton
tail is $O(r_*^4)\int_{9r_*}^{2R_{\rm ext}}r\,dr$.
The core replacement changes density only in $B_{7r_*}$,
where both approximate densities are bounded by one.
The selected initial perturbation is $O_{r_*}(\eta T^\delta)$.

By \eqref{eq:approx-subcritical} and the Hessian error, the initial
shifted energy is small. The unnormalized identities are
$$
\frac{d}{d\bar s}E(\bar f)=-\int_M\mathcal N(\bar f)^2d\mu_{\bar f},
\qquad
\frac{d}{d\bar s}\int_Me^{3\bar f}d\mu_0=-3E(\bar f).
$$
For fixed $r_*$, choose $T$ finally so that the shifted volume
stays between $1/2$ and $2$ and $\bar T$ is as small as required.
Normalization multiplies $\|R_{\bar f}^+\|_{L^\infty(M)}$ by
$\bar V^{1/3}$ and each curvature $L^{p_0}$ norm by
$\bar V^{1/3-1/p_0}$. Both factors are uniformly bounded.
The identity
$$
\|R_{\bar f}\|_{p_0}^{p_0}
=3^{p_0}\int_Me^{(3-p_0)\bar f}|\Delta_X\bar f|^{p_0}d\mu_0
$$
and the corresponding complex Hessian formulas for Ricci and
full Chern curvature prove their smallness.
The normalized lifetime is at most $2^{2/3}\bar T$.

The shifted normalization is exactly the original one:
$$
\bar V=e^{-3(c_\infty-C_{r_*})}V,\qquad
\bar f-\tfrac13\log\bar V=\widetilde f-\tfrac13\log V,\qquad
\bar V^{-2/3}d\bar s=V^{-2/3}ds.
$$
Thus it preserves the same profile and original trajectory.
Choose $r_*$ first, then all fixed spectral and core parameters
in the localization order, then $\eta$, and finally $T$.
Equations~\eqref{eq:exterior-negative} and
\eqref{eq:approx-subcritical}, followed by the actual error bounds,
give the prescribed $\varepsilon$ and $\tau_*$ simultaneously.
\end{proof}

\section{Curvature asymptotics and terminal geometry}\label{sec:constructed-geometry}

Normalize the localized solutions of Proposition~\ref{prop:smallness}
by \eqref{eq:normalize}. Denote their unnormalized and normalized
terminal times by $T$ and $\Theta$, respectively, and write
$t=t(s)$ for the time change. The profile $F$, its constants
$h,k,a,b$ and the model metric $\gamma_F$ are those defined in
Section~\ref{sec:shrinking}.

\begin{proposition}\label{prop:blowup-geometry}
For each of these solutions,
\begin{equation}\label{eq:exactrate}
\lim_{t\uparrow\Theta}\sqrt{\Theta-t}\,\|R_f(t)\|_\infty
=\sup_{\mathbb R^6}|R_F|\in(0,\infty),\qquad
\lim_{t\uparrow\Theta}\sqrt{\Theta-t}\,R_f(t,p)=-3h.
\end{equation}
For any $t_i\uparrow\Theta$, set $s_i=t^{-1}(t_i)$ and, in the
affine chart used in the construction,
$\Phi_i(y)=p+(T-s_i)^by$. Then
\begin{equation}\label{eq:normalizedlimit}
\Phi_i^*\!\left((\Theta-t_i)^{-1/2}
\gamma_f(\Theta+(\Theta-t_i)\theta)\right)
\longrightarrow
\tfrac12(-\theta)^a e^{F(|y|/(-\theta)^b)}g_{\mathrm E}
\end{equation}
in $C^\infty_{\rm loc}(\mathbb R^6\times(-\infty,0))$.
The limit is the unnormalized shrinking solution defined
in Section~\ref{sec:shrinking}; every slice is complete and nonflat.
The convergence holds for the full sequence with these affine maps.
\end{proposition}

\begin{proof}[Proof of Theorem~\ref{thm:construction},
Corollary~\ref{cor:smallness} and
Proposition~\ref{prop:blowup-geometry}]
Use Proposition~\ref{prop:smallness} and its actual localized solution
$\widetilde f(s)$, with terminal time $T$ before normalization.
Put $\sigma=T-s$ and
$V(s)=\int_Me^{3\widetilde f(s)}d\mu_0$.
The volume decreases and is bounded below, so $V(s)\to V_T>0$.
For the normalized time in \eqref{eq:normalize},
\begin{equation}\label{eq:remainingtime}
\frac{\Theta-t(s)}{T-s}\longrightarrow V_T^{-2/3}.
\end{equation}
Indeed the numerator is $\int_s^TV(u)^{-2/3}\,du$.

On a fixed rescaled ball,
$$
\sigma^{1/2}R_{\widetilde f}(s,\sigma^by)
=-3e^{-F-v(s,\sigma^by)}
\Delta_y(F+v(s,\sigma^by))\longrightarrow R_F.
$$
Here $a+2b=1/2$ and \eqref{eq:core-output} justify the convergence.
At $y=0$ its limit is $-3h$.
On $R\sigma^b\le|X|\le r_*$, the tail and
\eqref{eq:annular-output} give
$e^{F+v}\asymp |y|^k$ and
$|\nabla_y^{0,j}(F+v)|\le C|y|^{-j}$ for $j=1,2$.
Therefore
\begin{equation}\label{eq:scalar-tail}
\sigma^{1/2}|R_{\widetilde f}(s,X)|
+\sigma^{1/2}|\Ric^{\rm Ch}(\omega_{\widetilde f}(s))|
\le C|y|^{-k-2}.
\end{equation}
The fixed exterior contributes $o(1)$ after multiplication by
$\sqrt{\sigma}$. First take the time limit on $|y|\le R$, and then
let $R\to\infty$. This proves the global scalar supremum limit
and the analogous positive finite Ricci supremum limit.
Normalization multiplies curvature norms by $V^{1/3}$, while
\eqref{eq:remainingtime} gives
$\sqrt{(\Theta-t)/\sigma}\,V^{1/3}\to1$.
This proves both limits in \eqref{eq:exactrate}.

We justify smooth convergence in the original normalized time
without a general curvature compactness theorem.
Fix $0<\epsilon<A<\infty$. For every center $x$ set
$$
\ell_\sigma(x)=\sigma^{1/4}
e^{-\widetilde f(T-\sigma,x)/2},
\qquad
u_{\sigma,x}(\theta,y)=
\widetilde f(T+\sigma\theta,x+\ell_\sigma(x)y)
-\widetilde f(T-\sigma,x).
$$
For small $\sigma$, the affine maps embed $B_4$ at every center, and
on $B_4\times[-2A,-\epsilon/2]$ these functions have uniform
coefficient $C^2$ bounds and solve
$\partial_\theta u_{\sigma,x}=\mathcal N_y(u_{\sigma,x})$.
Here is the verification over all centers.
For $|x|\le R\sigma^b$, the ratio $\ell_\sigma(x)/\sigma^b$
is bounded above and below for fixed $R$.
The core expansion, after cancellation of $a\log\sigma$,
gives the bounds throughout a fixed enlarged rescaled ball and
the indicated negative time interval.
For an annular center of radius $r$, one has
$$
e^{\widetilde f(T-\sigma,x)}\asymp r^k,\qquad
\frac{\ell_\sigma(x)}r
\le C\left(\frac{\sigma^b}{r}\right)^{(k+2)/2}.
$$
Choose $R$ large depending on $A,\epsilon$.
When $r\ge R\sigma^b$, the spatial cylinder stays in a comparable
annulus at every indicated time, and the logarithmic time terms
cancel there. The annular zeroth-order and derivative bounds give
the required $C^2$ control.
For $x$ in the fixed exterior, $\ell_\sigma(x)\asymp\sigma^{1/4}$
and the cylinder stays in a slightly enlarged fixed exterior.
Its $C^2$ control follows from Proposition~\ref{prop:localized}.
These regions cover $M$, and
$\sup_x\ell_\sigma(x)\to0$ gives the uniform embedding.
Finally
$\ell_\sigma(x)^4=\sigma e^{-2\widetilde f(T-\sigma,x)}$,
so the spatial and constant-shift scaling factors in
\eqref{eq:N} cancel exactly. The chart is independent of $\theta$.

Lemma~\ref{lem:smoothing} on nested cylinders now gives all spatial
derivatives and the first time derivative. Repeated differentiation
of $\partial_\theta u=\mathcal N(u)$ expresses every higher time
derivative using finitely many spatial derivatives, so it gives
all mixed derivative bounds. In particular, for each compact
negative time interval and every integer $m\ge1$,
$$
\sup_{x,\theta}|\sigma^m\partial_s^m
\widetilde f(T+\sigma\theta,x)|\le C_m.
$$
Differentiating the volume density under its integral therefore gives
$|\partial_\theta^mV(T+\sigma\theta)|\le C_mV(0)$.
The functions $V(T+\sigma\theta)$ converge uniformly to $V_T$;
the next derivative bound and compactness show that this convergence
is $C^\infty$ on every compact negative interval.

Take arbitrary $t_i=t(s_i)\uparrow\Theta$, put
$\sigma_i=T-s_i$, and let $\Phi_i(y)=p+\sigma_i^by$.
Define the increasing time map
$$
\beta_i(\theta)=
\frac{t(T+\sigma_i\theta)-\Theta}{\Theta-t_i}.
$$
Equation~\eqref{eq:remainingtime} and the preceding volume convergence
give $\beta_i\to\operatorname{id}$ in $C^\infty_{\rm loc}(-\infty,0)$;
its derivative is
$\sigma_i(\Theta-t_i)^{-1}V(T+\sigma_i\theta)^{-2/3}\to1$.
Uniform convergence also holds up to the limiting zero endpoint
in the integral defining $\beta_i$.
The inverse maps converge smoothly on compact negative intervals,
by differentiating $\beta_i(\beta_i^{-1}(\theta))=\theta$.

The unnormalized pullbacks have the exact expression
\begin{align*}
\Phi_i^*(\sigma_i^{-1/2}\gamma_{\widetilde f(T+\sigma_i\theta)})
=\tfrac12(-\theta)^a\exp\left[
F(|y|/(-\theta)^b)+v(T+\sigma_i\theta,\sigma_i^by)\right]g_{\mathrm E}.
\end{align*}
The core bounds and the zero-forcing equation imply full smooth
convergence on compact space-time sets to the displayed model.
The normalized pullback at $\theta$ is the unnormalized pullback
at $\beta_i^{-1}(\theta)$ multiplied by
$$
\sqrt{\frac{\sigma_i}{\Theta-t_i}}\,
V(T+\sigma_i\beta_i^{-1}(\theta))^{-1/3}
\longrightarrow1
$$
smoothly. This proves \eqref{eq:normalizedlimit} for the full sequence.
The maps are holomorphic affine maps and preserve the fixed
complex structure.
Completeness of every limit slice follows from
$e^{F(r)/2}\asymp r^{k/2}$ at infinity: every escaping curve has
infinite radial lower length. Since $R_F(0)=-3h$, it is nonflat.

On each compact subset of $M\setminus\{p\}$ the unnormalized
solution has bounded amplitude and $C^2$ norm.
Lemma~\ref{lem:smoothing} and the local equation bound all spatial
derivatives and their first time derivatives at late times.
They are Cauchy in every spatial $C^j$ norm. The positive volume
limit preserves a positive smooth spatial limit after normalization.
At $p$, $\widetilde f(s,p)=a\log\sigma+o(1)\to-\infty$,
so no positive smooth extension exists. Uniqueness identifies
$\Theta$ as the maximal normalized time and $p$ as its only
singular point. Since $h>0$, the two limits in
\eqref{eq:exactrate} give the curvature assertions of
Theorem~\ref{thm:construction}.
Proposition~\ref{prop:smallness} supplies the prescribed smallness
and lifetime in Corollary~\ref{cor:smallness}.
\end{proof}

\begin{corollary}\label{cor:terminal}
For every solution in Proposition~\ref{prop:blowup-geometry}, the completion
of $(M\setminus\{p\},\gamma_\Theta)$ adds exactly one point.
The compact spaces $(M,\gamma_f(t))$ converge to this completion
in Gromov--Hausdorff distance. The added point has the unique
nonflat metric tangent cone
$$
d\rho^2+\frac{(k+2)^2}{4}\rho^2g_{S^5}.
$$
Writing $R_{\gamma_\Theta}$ for its Chern scalar curvature, one has
near the added point
$$
\rho^2R_{\gamma_\Theta}
=-\frac{24k}{(k+2)^2}+O(\rho^{4\delta})
$$
for some $\delta>0$, with the differentiated control proved in
Section~\ref{sec:constructed-geometry}.
\end{corollary}

\begin{proof}[Proof of Corollary~\ref{cor:terminal}]
The physical annular bounds pass to the punctured terminal limit,
and give
$$
\gamma_\Theta=\tfrac12V_T^{-1/3}e^{c_\infty}
r^ke^{v_T(X)}(dr^2+r^2g_{S^5}),\qquad
|\nabla^{0,j}v_T|\le Cr^{\delta/b-j}\quad(0\le j\le2).
$$
Set
$$
\rho=\frac{\sqrt2V_T^{-1/6}e^{c_\infty/2}}{k+2}r^{(k+2)/2}.
$$
Direct substitution yields
$$
\gamma_\Theta=e^{v_T}
\left(d\rho^2+\frac{(k+2)^2}{4}\rho^2g_{S^5}\right),
\qquad v_T=O(\rho^{4\delta}),
$$
with the corresponding two weighted derivatives.
The distance to $r=0$ is finite and the diameter of the completed
$B_r$ is at most $Cr^{(k+2)/2}$.
A separating closed annulus is uniformly positive and keeps all
other points separated. The completion therefore adds exactly one
point.

Under $\rho=\lambda R$ and metric multiplication by $\lambda^{-2}$,
the conformal error tends to zero on each closed annulus in $R>0$.
The inner ball has rescaled diameter at most $CR$, uniformly in
$\lambda$, so its diameter tends to zero with $R$.
This proves convergence including the vertex to the stated cone
for all $\lambda\downarrow0$, hence uniqueness.
Its tangential sectional curvature is
$[4-(k+2)^2]/[(k+2)^2\rho^2]<0$, so the metric cone is nonflat.
The scalar formula $R=-3e^{-f}\Delta_Xf$ applied to the terminal
expression gives
$$
\rho^2R_{\gamma_\Theta}
=-\frac{24k}{(k+2)^2}+O(\rho^{4\delta}).
$$

For the compact convergence, the core, tail and volume bounds imply
$$
e^{f(t(s),X)}\le C(\sigma^a+|X|^k).
$$
Whenever $\sigma^b\le r$, radial paths through $p$ give preterminal
diameter at most $Cr^{1+k/2}$ for $B_r$.
The terminal completion has the same bound.
Collapse $B_r$ to one point in both spaces. This changes each
distance by at most that diameter: a path through the collapsed
point can be lifted by connecting its entrance and exit inside
the ball. On the common compact complement, smooth metric
convergence gives convergence of the length distances with the
boundary collapsed. These collapsed compact spaces converge in
GH distance. First let $t\uparrow\Theta$, then $r\downarrow0$;
both diameter errors vanish. This proves the corollary.
\end{proof}

\begin{proposition}\label{prop:compatibility}
The solutions in Proposition~\ref{prop:blowup-geometry} satisfy
\eqref{eq:scalar-log}, and \eqref{eq:scalar-gauge} for every $h$ in
Remark~\ref{rem:scalar-gauge}, with $T=\Theta$. They also satisfy
\eqref{eq:ricci-budget}. More precisely,
\begin{equation}\label{eq:ricci-L9-limit}
\lim_{t\uparrow\Theta}(\Theta-t)^{1/3}
\|\Ric^{\rm Ch}(\omega_f(t))\|_{L^9(M,d\mu_f(t))}
=\|\Ric^{\rm Ch}(\gamma_F)\|_{L^9(\mathbb R^6,d\mu_F)}
\in(0,\infty),
\end{equation}
where $\gamma_F=\tfrac12e^Fg_{\mathrm E}$ and $d\mu_F=e^{3F}dy/8$.
Consequently their terminal dyadic budgets converge to
$\|\Ric^{\rm Ch}(\gamma_F)\|_{L^9(d\mu_F)}^3\log2>0$.
\end{proposition}

\begin{proof}
Use the unnormalized construction time $s$, put $\sigma=T-s$,
and write $X=\sigma^by$.
The core and annular estimates, together with the bounded
normalization factors, give on the fixed construction chart
$$
e^{f(t(s),X)}\le C(\sigma^b+|X|)^k,\qquad
|R_f(t(s),X)|\le C(\sigma^b+|X|)^{-(k+2)}.
$$
The exterior has uniformly bounded curvature. For a sufficiently
large level $\lambda$, a nonempty superlevel set is contained in
$|X|<(C/\lambda)^{1/(k+2)}$, and nonemptiness also forces
$\sigma^b<(C/\lambda)^{1/(k+2)}$. Hence
$$
\mu_{f(t)}\{|R_f(t)|>\lambda\}
\le C\lambda^{-(6+3k)/(k+2)}=C\lambda^{-3}.
$$
Lower levels are controlled by unit volume, and earlier compact
time intervals are smooth. Thus the weak $L^3$ curvature quasi-norm
is uniformly bounded. The Type-I bound in
Proposition~\ref{prop:blowup-geometry} controls
$\sqrt{\Theta-t}\,\|R_f(t)\|_\infty$.
For $\lambda\le\|R_f(t)\|_\infty$, local boundedness and monotonicity
of $h$ bound $h(\sqrt{\Theta-t}\lambda)$ by a fixed constant;
larger levels have empty superlevel sets. This proves
\eqref{eq:scalar-gauge}, and therefore \eqref{eq:scalar-log}.

The actual core convergence gives
$$
\sigma^{1/2}|\Ric^{\rm Ch}(\omega_{\widetilde f})|
(s,\sigma^by)\longrightarrow|\Ric^{\rm Ch}(\gamma_F)|(y),
\qquad
\sigma^{-3/2}\Phi_\sigma^*d\mu_{\widetilde f}
\longrightarrow d\mu_F
$$
on every fixed rescaled ball, where $\Phi_\sigma(y)=\sigma^by$.
The volume exponent is $3a+6b=3/2$.
The annular Hessian estimate in \eqref{eq:annular-output} and
$e^{F(r)}\asymp r^k$ give, for every sufficiently large fixed $R$,
\begin{align*}
&\sigma^3\int_{R\sigma^b\le|X|\le r_*}
|\Ric^{\rm Ch}(\omega_{\widetilde f})|^9d\mu_{\widetilde f}\\
&\qquad\le C\int_R^\infty r^{5+3k-9(k+2)}\,dr
\le CR^{-6(k+2)}.
\end{align*}
The fixed exterior has uniformly bounded curvature and volume,
so its integral multiplied by $\sigma^3$ tends to zero.
First taking $s\uparrow T$ on each fixed core and then
$R\uparrow\infty$ yields
$$
\sigma^3\int_M|\Ric^{\rm Ch}(\omega_{\widetilde f})|^9
\,d\mu_{\widetilde f}
\longrightarrow
\int_{\mathbb R^6}|\Ric^{\rm Ch}(\gamma_F)|^9\,d\mu_F.
$$
The limit is finite by the same tail estimate and positive because
$R_F(0)\ne0$.

Under \eqref{eq:normalize}, the ninth-power integral is multiplied
by $V^2$. Equation~\eqref{eq:remainingtime} gives
$((\Theta-t)/\sigma)^3V^2\to1$, proving
\eqref{eq:ricci-L9-limit}. Cubing that limit gives the curvature
integrand as $(\Theta-t)^{-1}$ times a function with a positive
finite limit. Integration over $[\Theta-2s,\Theta-s]$ supplies the
stated limiting budget. All earlier compact time intervals are
smooth, so \eqref{eq:ricci-budget} holds.
\end{proof}

\begin{proposition}\label{prop:critical}
For every fixed solution in Proposition~\ref{prop:blowup-geometry},
$$
\int_M|R_f(t)|^3d\mu_f
=\frac{108\pi^3k^3}{k+2}\log\frac1{\Theta-t}+O(1).
$$
\end{proposition}

\begin{proof}
At the critical exponent,
$\int|R_f|^3d\mu_f=27\int|\Delta_Xf|^3d\mu_0$,
unchanged by constant normalization.
On $R\sigma^b\le r\le r_*$,
$$
r^2\Delta_X\widetilde f
=4k+O(\sigma r^{-2(k+2)})+O(r^{\delta/b}).
$$
The cubic difference from $(4k)^3$ is bounded by the sum of
the displayed errors. Its radial integral against $dr/r$ is
uniformly bounded: the first term integrates to at most a
constant times $R^{-2(k+2)}$, and the second to a fixed constant.
As $|S^5|=\pi^3$ and $d\mu_0=dX/8$, the annulus contributes
$$
\frac{27}{8}\pi^3(4k)^3\log\frac{r_*}{R\sigma^b}+O(1).
$$
The core contribution is bounded by rescaling $X=\sigma^by$,
and the exterior contribution is bounded by its fixed $C^2$ norm.
Since $b=1/[2(k+2)]$ and
$\log(1/\sigma)=\log(1/(\Theta-t))+O(1)$ by
\eqref{eq:remainingtime}, the stated coefficient follows.
\end{proof}

\begin{proposition}\label{prop:terminal-energy}
For the same solution let $E_\Theta=\lim_{t\uparrow\Theta}E(t)$.
Then
$$
\int_{M\setminus\{p\}}d\mu_{\gamma_\Theta}=1,\qquad
E_\Theta=\int_{M\setminus\{p\}}R_{\gamma_\Theta}^2d\mu_{\gamma_\Theta}
\in(0,\infty).
$$
The volume and curvature-square measures converge weakly, in the
fixed background identification, to these punctured measures
extended by zero at $p$. Moreover
\begin{align}
\lim_{t\uparrow\Theta}\sqrt{\Theta-t}\,[-E'(t)]
&=\int_{\mathbb R^6}(-a+b|y|F'(|y|))^2d\mu_F,\label{eq:dissipation-limit}\\
\lim_{t\uparrow\Theta}\frac{E(t)-E_\Theta}{\sqrt{\Theta-t}}
&=2\int_{\mathbb R^6}(-a+b|y|F'(|y|))^2d\mu_F.\label{eq:energy-defect}
\end{align}
The profile integral is finite and strictly positive.
\end{proposition}

\begin{proof}
For $r<r_0$ and all sufficiently late times with $R_0\sigma^b\le r$,
the core and physical annulus bounds give
\begin{equation}\label{eq:smallballs}
\mu_f(B_r)\le Cr^{3k+6},\qquad
\int_{B_r}R_f^2d\mu_f\le Cr^{k+2}.
\end{equation}
Indeed the core volume is at most $C\sigma^{3a+6b}=C\sigma^{3/2}$
and its energy at most $C\sigma^{1/2}$.
These are bounded by the respective powers of $r$.
The annular radial integrands are at most $Cr^{3k+5}$ and
$Cr^{k+1}$. The bounded normalization factors preserve the bounds.
Smooth convergence outside $B_r$, followed by $r\downarrow0$,
proves the total-volume and energy identities and weak convergence.
The terminal scalar asymptotic is nonzero on a fixed small annulus,
so $E_\Theta>0$. The restriction $R_0\sigma^b\le r$ in
\eqref{eq:smallballs} is sufficient for this limiting argument.

The profile equation gives
$\mathcal N(F)=-a+brF'$.
It is nonzero at the origin and has order $r^{-2(k+2)}$ at infinity.
Against $d\mu_F\asymp r^{3k+5}dr$ its square has integrable tail
$r^{-k-3}dr$. Thus its squared integral is finite and positive.

On the actual core, the PDE, rather than time differentiation of
an error estimate, gives
$$
\sigma\widetilde f_s(s,\sigma^by)
=\mathcal N(F+v(s,\sigma^b\cdot))(y).
$$
The pulled-back volume is $\sigma^{3/2}e^{3(F+v)}dy/8$.
Consequently, on every fixed core,
$\sqrt{\sigma}\int\widetilde f_s^2d\mu_{\widetilde f}$
converges to the corresponding profile integral.
On the annulus, $|\widetilde f_s|\le Cr^{-2k-4}$ follows directly
from \eqref{eq:N} and four annular derivatives, and hence
$$
\sqrt{\sigma}\int_{R\sigma^b\le|X|\le r_*/2}
\widetilde f_s^2d\mu_{\widetilde f}
\le C\sqrt{\sigma}\int_{R\sigma^b}^{r_*/2}r^{-k-3}\,dr
\le CR^{-k-2}.
$$
The exterior contribution tends to zero. This proves the global
unnormalized dissipation limit.

Normalization gives
$f_t=V^{2/3}\widetilde f_s+E_f$ and
$\int\widetilde f_s\,d\mu_{\widetilde f}=-E(\widetilde f)$.
Expanding the square, including its cross term, yields
$$
\int_Mf_t^2d\mu_f
=V^{1/3}\int_M\widetilde f_s^2d\mu_{\widetilde f}-E_f^2.
$$
The bounded last term disappears after multiplication by
$\sqrt{\Theta-t}$, and \eqref{eq:remainingtime} cancels the
remaining volume factor. Equation~\eqref{eq:dissipation}
proves \eqref{eq:dissipation-limit}. Integrating its two-sided
asymptotic estimate from $t$ to $\Theta$ proves
\eqref{eq:energy-defect}.
The finite terminal energy here is distinct from the full
model integral $\int_{\mathbb R^6}R_F^2d\mu_F$, whose tail
$r^{k+1}dr$ is not integrable.
\end{proof}

\appendix

\section{Proof of the radial profile construction}\label{sec:profile}

We use the real Euclidean Laplacian on $\mathbb R^6$.
Recall that $\kappa=\sqrt5-1$ and put
$$
\Omega=\left(\frac{23\sqrt5-43}{2}\right)^{1/2}.
$$
The parameter $k$ varies near $\kappa$ in this appendix.

\begin{proof}[Proof of Lemma~\ref{lem:profile}]
\medskip\noindent\textbf{The radial equation and the regular origin.}
For $k$ close to $\kappa$, put $b=1/[2(k+2)]$ and $a=bk$. A radial profile must satisfy
$$
-18e^{-3F}\Delta(e^F\Delta F)-9e^{-2F}(\Delta F)^2=-a+brF_r.
$$
Let $H=e^F\Delta F$ and $p=rF_r$. This equation is equivalent to
$$
\Delta H=-\frac12H^2e^{-F}+\frac1{18}(a-bp)e^{3F}.
$$
The regular solution with $F(0)=0$ and $H(0)=h$ is determined by
\begin{align*}
p(r)&=r^{-4}\int_0^r t^5H(t)e^{-F(t)}\,dt,\\
F(r)&=\int_0^r\frac{p(t)}t\,dt,\\
H(r)&=h+\int_0^r t^{-5}\int_0^t
s^5\left[-\frac12H(s)^2e^{-F(s)}
+\frac1{18}(a-bp(s))e^{3F(s)}\right]\,ds\,dt.
\end{align*}
These integral equations also establish regularity at $r=0$. Indeed, in the variable $x=r^2$, if $S(x)$ denotes the expression in square brackets, they give
$$
H_x(x)=\frac14\int_0^1t^2S(xt)\,dt,\qquad
F_x(x)=\frac14\int_0^1t^2H(xt)e^{-F(xt)}\,dt,\qquad
p=2xF_x.
$$
In the right-hand side defining $S$, the last expression for $p$ can be replaced by the preceding integral for $F_x$. The integrated equations are therefore a contraction on a small $x$-interval in a ball of continuous pairs $(F,H)$ around $(0,h)$. The integral operators are smooth functions of their arguments. Starting with continuous $F,H$, the displayed formulas first give one derivative and then, by differentiation, every derivative. Thus $F$ is a smooth function of $r^2$. The same contraction and differentiation give smooth dependence on $h$ and $k$. In particular,
$$
p(r)=\frac h6r^2+O(r^4).
$$

For $r>0$, introduce
$$
s=\log r,\qquad p=F_s,\qquad z=F_{ss}+4F_s,\qquad
w=z_s,\qquad A=e^{2F+4s}.
$$
Direct differentiation gives the exact system
\begin{align}
p_s&=z-4p,\qquad z_s=w,\notag\\
w_s&=-2pw-(p^2-4p-4)z-\frac32z^2
+\frac b{18}(k-p)A,\notag\\
A_s&=2(p+2)A. \label{eq:profile-system}
\end{align}
The regular origin has
$$
p=\frac h6r^2+O(r^4),\quad z=hr^2+O(r^4),\quad
w=2hr^2+O(r^4),\quad A=r^4(1+O(r^2)).
$$

\medskip\noindent\textbf{A stationary inner solution.}
We first omit the term containing $A$ in the first three equations of
\eqref{eq:profile-system}. Equivalently, we solve the stationary equation
$$
-18e^{-3F}\Delta(e^F\Delta F)-9e^{-2F}(\Delta F)^2=0
$$
with regular data $F(0)=0$, $\Delta F(0)=1$. Its local existence follows from the same integral equations, with $S=-H^2e^{-F}/2$.

Direct substitution gives
$$
\frac d{ds}\left[e^{F+2s}
\bigl(z^2+(8-8p-2p^2)z-2(p+2)w\bigr)\right]=0.
$$
Thus the displayed quantity is constant. Its limit at the regular origin is zero. Hence $z^2+(8-8p-2p^2)z-2(p+2)w=0$, and the stationary orbit satisfies
$$
p_s=z-4p,\qquad
z_s=\frac{z(z+8-8p-2p^2)}{2(p+2)}.
$$
The origin expansion has $p>0$ and $z>0$. To prove global existence and determine the limit, put
$$
v=p+2,\qquad Y=\frac z{p+2}.
$$
The system becomes
$$
v_s=8+v(Y-4),\qquad Y_s=Y(4-v-Y/2).
$$
The region $v>2$, $Y>0$ is forward invariant: on $v=2$ the first component is $2Y>0$, and the second equation preserves $Y>0$. Its positive equilibrium is
$$
v_*=(\kappa+2),\qquad Y_*=6-2\sqrt5.
$$
The function
$$
v-v_*-v_*\log(v/v_*)
+Y-Y_*-Y_*\log(Y/Y_*)
$$
is proper on the positive quadrant, and a direct calculation gives
$$
\frac d{ds}\left[v-v_*-v_*\log(v/v_*)+Y-Y_*-Y_*\log(Y/Y_*)\right]
=-\frac{8(v-v_*)^2}{vv_*}-\frac12(Y-Y_*)^2.
$$
Starting at any finite $s$, the orbit therefore remains in a compact subset of the positive quadrant and exists for all later $s$. The two nonnegative terms on the right are integrable. Their derivatives are bounded on this compact set, so they tend to zero. It follows that
$$
p\longrightarrow\kappa,\qquad z\longrightarrow4\kappa.
$$
We denote this stationary profile by $F_*$. Since $p>0$, it is strictly increasing.

The eigenvalues of the planar system at this equilibrium are
$-\frac{\kappa+2}{2}\pm i\Omega$. The remaining eigenvalue of the full three-dimensional stationary system is $-(\kappa+2)$. The first two eigenvalues follow from the characteristic polynomial
$$
\lambda^2+(\kappa+2)\lambda+12\sqrt5-20.
$$
In real linear coordinates for the planar focus, write the orbit as $(R,\theta)$. For sufficiently large $s$,
$$
R_s=-\frac{\kappa+2}{2} R+O(R^2),\qquad \theta_s=\Omega+O(R).
$$
The orbit is not the equilibrium, by its regular-origin expansion. Hence $R>0$. First these equations give exponential decay with any slightly smaller rate than $\frac{\kappa+2}{2}$. Then $R$ is integrable, and integration of
$(\log R)_s=-\frac{\kappa+2}{2}+O(R)$ shows that $e^{\frac{\kappa+2}{2} s}R$ has a finite strictly positive limit. Integration of the angular equation gives a finite limit of $\theta-\Omega s$.

Writing
$$
U_*(s)=(p-\kappa,z-4\kappa,w),
$$
we obtain a nonzero complex coefficient $c$ and an eigenvector $v_+$ such that
$$
U_*(s)=e^{-\frac{\kappa+2}{2} s}\operatorname{Re}(ce^{i\Omega s}v_+)
+O((1+s)e^{-(\kappa+2) s}).
$$
For completeness, the stated error follows by variation of constants after the leading oscillatory coefficients have been subtracted. The nonlinear forcing is $O(e^{-(\kappa+2) s})$. Its convolution in the real eigendirection has size $O((1+s)e^{-(\kappa+2)s})$, because the forcing and the real eigendirection have the same decay exponent. The oscillatory eigendirections give the same or a better bound. Thus no nonresonance assumption is being made. Integration of the first component gives
$$
F_*(r)=\kappa\log r+c_{\rm in}+O(r^{-\frac{\kappa+2}{2}}).
$$
For arbitrary $h>0$, the stationary profile with $\Delta F(0)=h$ is
$F_*(\sqrt h\,r)$.

\medskip\noindent\textbf{The tail condition at infinity.}
We next construct the solutions of \eqref{eq:profile-system} which have
$p\to k$ at infinity. Eliminating $z$ and $w$ gives
$$
p_{sss}+2(p+2)p_{ss}+(p^2+16p-4)p_s+\frac32p_s^2
+4p(p^2+2p-4)+\frac b{18}(p-k)A=0.
$$
As long as $p>-2$, the variable $A$ increases. Put $q=A^{1/3}$. Since
$D_s=\frac23(p+2)qD_q$, the preceding equation becomes
\begin{align}
0={}&p_{qqq}+\frac6q p_{qq}
+\left(4+\frac{9(p^2+16p-4)}{4(p+2)^2}\right)\frac{p_q}{q^2}
+\frac4{p+2}p_qp_{qq}\notag\\
&+\frac{37}{4(p+2)q}p_q^2+\frac{p_q^3}{(p+2)^2}
+\frac{27p(p^2+2p-4)}{2(p+2)^3q^3}
+\frac{3b(p-k)}{16(p+2)^3}. \label{eq:profile-tail-equation}
\end{align}
This identity follows by expanding the three applications of
$\frac23(p+2)qD_q$ in the scalar equation; in particular, it retains the complete nonlinear equation.

Set $X=(p-k,p_q,p_{qq})$. Equation \eqref{eq:profile-tail-equation} has the form
$$
X_q=B_kX+\mathcal F(q,k,X),\qquad
B_k=\begin{pmatrix}0&1&0\\0&0&1\\-\mu_k^3&0&0\end{pmatrix},
\qquad \mu_k^3=\frac3{32(k+2)^4}.
$$
For $k$ in a fixed small interval around $\kappa$, $q\geq q_0$, and $X$ in a fixed small ball, its explicit terms give
\begin{align*}
|\mathcal F(q,k,X)|
&\leq C\big(q^{-1}|X|+|X|^2+|k-\kappa|q^{-3}\big),\\
|\mathcal F(q,k,X)-\mathcal F(q,k,\widetilde X)|
&\leq C\big(q^{-1}+|X|+|\widetilde X|\big)|X-\widetilde X|.
\end{align*}
Here the constant term vanishes at $k=\kappa$ because
$\kappa^2+2\kappa-4=0$. The eigenvalues of $B_k$ are
$-\mu_k$ and $\mu_k/2\pm i\sqrt3\mu_k/2$. Let $\Pi_s,\Pi_u$ be its real stable and unstable projections, and let $e_s$ be a unit stable eigenvector.

For a small real parameter $\sigma$ solve
\begin{align*}
X(q)={}&e^{B_k(q-q_0)}\sigma e_s
+\int_{q_0}^q e^{B_k(q-t)}\Pi_s\mathcal F(t,k,X(t))\,dt\\
&-\int_q^\infty e^{B_k(q-t)}\Pi_u\mathcal F(t,k,X(t))\,dt.
\end{align*}
In the norm $\|X\|=\sup_{q\geq q_0}q^3|X(q)|$, the right-hand side defines an operator $\mathcal T$ satisfying
\begin{align*}
\|\mathcal T X\|
&\leq C(q_0^3|\sigma|+|k-\kappa|)
+C(q_0^{-1}\|X\|+q_0^{-3}\|X\|^2),\\
\|\mathcal T X-\mathcal T\widetilde X\|
&\leq C\big(q_0^{-1}+q_0^{-3}(\|X\|+\|\widetilde X\|)\big)
\|X-\widetilde X\|.
\end{align*}
The weighted convolution estimates can be seen by splitting the stable integral at $t=q/2$. Above that point, the polynomial weights are comparable; below it, the exponential kernel absorbs every polynomial factor. In the unstable integral $t\geq q$, so the weight ratio is at most one. Choose $q_0$ large and then $q_0^3|\sigma|+|k-\kappa|$ small. The contraction theorem gives a unique tail satisfying $X=O(q^{-3})$. Differentiating the same contraction equation gives smooth dependence on $k$ and $\sigma$.

Moreover, $\Pi_sX(q_0)=\sigma e_s$ exactly. Thus these tails form an immersed one-dimensional curve, and they comprise all sufficiently small tails with $X=O(q^{-3})$. Flowing backward over a fixed interval of positive $q$ gives a smooth curve of codimension two at $A=1$, depending smoothly on $k$. At $k=\kappa$ and $\sigma=0$, it is the exact cone
$$
p=\kappa,\qquad z=4\kappa,\qquad w=0.
$$

\medskip\noindent\textbf{Two linear connection facts.}
We verify the two transversality assertions needed to join the regular origin to the tail. Linearizing the $p$ equation at the exact cone, with $k=\kappa$, gives
$$
\left[(D_s+(\kappa+2))(D_s^2+(\kappa+2)D_s+12\sqrt5-20)
+\frac{A}{36(\kappa+2)}\right]g=0.
$$
Put
$$
x=\frac{A}{288(\kappa+2)^4},\qquad \theta=x\frac d{dx},\qquad
P(\theta)=(\theta+\tfrac12)(\theta+\tfrac14-i\frac{\Omega}{2(\kappa+2)})(\theta+\tfrac14+i\frac{\Omega}{2(\kappa+2)}).
$$
The equation is $(P(\theta)+x)g=0$. Its exponents at zero are
$-1/2$ and $-1/4\pm i\frac{\Omega}{2(\kappa+2)}$.

Consider the real-valued function
$$
D(x)=x^{-1/2}\int_0^\infty\!\!\int_0^\infty
t^{-3/4+i\frac{\Omega}{2(\kappa+2)}}u^{-3/4-i\frac{\Omega}{2(\kappa+2)}}
\exp\!\left(-t-u-\frac{x}{tu}\right)dt\,du.
$$
Its reality follows by conjugation and exchange of $t,u$. The integral and all its $\theta$ derivatives converge absolutely for $x>0$. Integration by parts in $t,u$ gives $(P(\theta)+x)D=0$. Here is the elementary rule behind this identity. In an integral with prefactor $x^{\beta_3}$ and powers
$t^{\beta_1-\beta_3-1}u^{\beta_2-\beta_3-1}$, successive application of
$\theta-\beta_1$, $\theta-\beta_2$, and $\theta-\beta_3$ inserts the multipliers $-t$, $-u$, and $-x/(tu)$, respectively. Their product is $-x$. For the displayed integral take
$(\beta_1,\beta_2,\beta_3)=(-1/4+i\frac{\Omega}{2(\kappa+2)},-1/4-i\frac{\Omega}{2(\kappa+2)},-1/2)$.
All integration-by-parts boundary terms vanish exponentially.

Dominated convergence gives
$$
\lim_{x\downarrow0}x^{1/2}D(x)
=\left|\Gamma\!\left(\tfrac14+i\frac{\Omega}{2(\kappa+2)}\right)\right|^2>0;
$$
the absolute dominating function is $t^{-3/4}u^{-3/4}e^{-t-u}$.
The nonvanishing can also be read directly from Euler's product:
$$
|\Gamma(c+i\frac{\Omega}{2(\kappa+2)})|^2=\Gamma(c)^2
\prod_{n=0}^\infty\left(1+\frac{\frac{\Omega}{2(\kappa+2)}^2}{(n+c)^2}\right)^{-1}>0,
\qquad c>0.
$$
The product is positive because the sum of its logarithms converges.
To estimate the other endpoint, use $t+u+x/(tu)\ge3x^{1/3}$.
Every $\theta^jD$ has the same integral with an additional polynomial in $x/(tu)$. Splitting the exponential into two equal factors and using the boundedness of $(x/(tu))^j e^{-x/(2tu)}$ yields
$$
|\theta^jD(x)|\le C_jx^{-1/2}e^{-3x^{1/3}/2},\qquad x\ge1.
$$
Thus $D$ is a nonzero member of the decaying linear tail space, which is one-dimensional by the linearized tail contraction. It spans that line. Every pure oscillatory Frobenius solution satisfies $x^{1/2}g(x)\to0$ at zero, whereas $D$ does not. Therefore the two-dimensional real space of oscillatory incoming modes projects isomorphically onto the quotient by the tail line. This is the first connection assertion.

For the second assertion, variation of $k$ by one gives the forced equation
$$
(P(\theta)+x)g=x.
$$
We show directly that no solution of this equation has $g,\theta g,\theta^2g$ bounded at both endpoints. Set
$$
\alpha_1=\frac12,\qquad \alpha_2=\frac14+i\frac{\Omega}{2(\kappa+2)},\qquad
\alpha_3=\frac14-i\frac{\Omega}{2(\kappa+2)},\qquad \zeta=e^{i\pi/3},
$$
and define
$$
\psi(x)=\zeta x^{\alpha_3}
\int_0^\infty\!\!\int_0^\infty
t^{\alpha_1-\alpha_3-1}u^{\alpha_2-\alpha_3-1}
\exp\!\left[-\zeta\left(t+u+\frac{x}{tu}\right)\right]dt\,du.
$$
Positive real powers use real logarithms, and $\zeta^z$ uses $\log\zeta=i\pi/3$. Since $\operatorname{Re}\zeta=1/2$, this integral and its derivatives converge absolutely. The same integration-by-parts rule now inserts $-\zeta t$, $-\zeta u$, and $-\zeta x/(tu)$. Since $-\zeta^3=1$, it gives
$$
\prod_{j=1}^3(\theta-\alpha_j)\psi=x\psi,
\qquad (P(-\theta)+x)\psi=0.
$$
The last operator is the formal transpose for the bilinear pairing with measure $dx/x$; no complex conjugation is used.

For $\operatorname{Re}s>-1/4$, integrating first in $x$ and then in $t,u$ gives
\begin{equation}\label{eq:adjoint-mellin-moment}
\int_0^\infty x^{s-1}\psi(x)\,dx
=\zeta^{-3s}\prod_{j=1}^3\Gamma(s+\alpha_j)
=e^{-i\pi s}\prod_{j=1}^3\Gamma(s+\alpha_j).
\end{equation}
Indeed, the first integral contributes
$\Gamma(s+\alpha_3)(tu/\zeta)^{s+\alpha_3}$, after which the powers of $t,u$ are $s+\alpha_1-1,s+\alpha_2-1$. Absolute Fubini is valid: replacing $\zeta$ in the exponential by its positive real part gives three Euler integrals with positive real arguments. The same absolute integrability holds with any $\theta^j\psi$ in place of $\psi$. Such differentiation inserts only a polynomial in $x/(tu)$, whose added powers cancel from the $t,u$ exponents after the $x$ integration.

Taking $s=0$ shows that all $\theta^j\psi$ belong to $L^1(dx/x)$. As functions of $\log x$, they therefore belong to $W^{1,1}(\mathbb R)$ and tend to zero at both endpoints. Taking $s=1$ in \eqref{eq:adjoint-mellin-moment} gives
$$
\int_0^\infty\psi(x)\,dx
=-\Gamma\!\left(\tfrac32\right)
\left|\Gamma\!\left(\tfrac54+i\frac{\Omega}{2(\kappa+2)}\right)\right|^2\ne0.
$$
If $g,\theta g,\theta^2g$ were bounded at both endpoints, three integrations by parts, with the boundary terms vanishing by the preceding decay, would give
$$
0=\int_0^\infty g(P(-\theta)+x)\psi\,\frac{dx}{x}
=\int_0^\infty\psi(P(\theta)+x)g\,\frac{dx}{x}
=\int_0^\infty\psi(x)\,dx,
$$
a contradiction. This proves the second connection assertion using only the displayed integrals.

\medskip\noindent\textbf{A curve leaving the cone.}
For the nonlinear matching, write
$$
U=(p-\kappa,z-4\kappa,w).
$$
Equation \eqref{eq:profile-system} becomes
$$
U_s=\mathsf J U+\mathcal R_0(U)
+\frac{b(k)}{18}(k-\kappa-U_1)Ae_3,\qquad
A_s=2((\kappa+2)+U_1)A,
$$
where $\mathsf J$ is the Jacobian with respect to $(p,z,w)$ of the first three right-hand sides of \eqref{eq:profile-system} at $(p,z,w,A)=(\kappa,4\kappa,0,0)$ with $k=\kappa$, and $\mathcal R_0$ is the polynomial remainder, with $\mathcal R_0(U)=O(|U|^2)$. The matrix $\mathsf J$ has the three eigenvalues $-(\kappa+2),-\frac{\kappa+2}{2}\pm i\Omega$.
In particular, $\mathsf J$ does not depend on $k$.

For $k$ near $\kappa$ there is a curve $U=\mathcal Q_k(A)$ with
$\mathcal Q_k(0)=0$ satisfying
$$
2((\kappa+2)+(\mathcal Q_k)_1)A\mathcal Q_k'
=\mathsf J\mathcal Q_k+\mathcal R_0(\mathcal Q_k)
+\frac{b(k)}{18}(k-\kappa-(\mathcal Q_k)_1)Ae_3.
$$
Here is a direct construction. Write its equation as
$$
\mathcal Q_k'=\frac{\mathsf J}{2(\kappa+2)A}\mathcal Q_k
+\frac{\mathcal G(A,k,\mathcal Q_k)}A,
$$
where
$$
|\mathcal G(A,k,Q)|\leq C(|Q|^2+A|Q|+|k-\kappa|A).
$$
Solve the integral equation
$$
\mathcal Q_k(A)=\int_0^A
\exp\!\left(\frac{\mathsf J}{2(\kappa+2)}\log\frac At\right)
\mathcal G(t,k,\mathcal Q_k(t))\,\frac{dt}t.
$$
The matrix exponential is bounded by $C(A/t)^{-1/4}$. In the norm
$\sup_{0<A\leq A_0}|\mathcal Q_k(A)|/A$, this is a contraction on a ball of radius $C|k-\kappa|$, after choosing $A_0$ small and then $k$ close to $\kappa$. It follows that
$$
\mathcal Q_k(A)=O(|k-\kappa|A),\qquad \mathcal Q_\kappa=0.
$$
Differentiating the integral equation proves smooth parameter dependence; the coefficients are analytic in $A$ by the same recursion. The curve extends to $0\leq A\leq1$ by ordinary differential equations, while remaining small for $k$ sufficiently close to $\kappa$.

Subtract this curve and put $V=U-\mathcal Q_k(A)$. Its equation has no additive error:
$$
V_s=\mathsf J V+\mathcal R_0(V)+\mathcal R(A,k,V),\qquad
|\mathcal R(A,k,V)|\leq CA|V|
$$
in a fixed small neighborhood. This estimate includes the term obtained by differentiating $\mathcal Q_k(A)$ along the actual equation for $A$.

\medskip\noindent\textbf{Passage from the regular origin.}
Return to the genuine regular-origin profile with $\Delta F(0)=h$, and scale
$$
\rho=\sqrt h\,r,\qquad F(r)=G_h(\rho).
$$
Its data are $G_h(0)=0$, $\Delta_\rho G_h(0)=1$.
In the regular integral equations, after also scaling $H$ by $h$, the term
$(a-bp)e^{3G_h}/18$ is multiplied by $h^{-2}$, while all stationary terms remain unchanged. Thus the contraction at the origin and ordinary differential equations on a fixed later interval give
$$
G_h-F_*=O(h^{-2})
$$
in every fixed finite differentiability norm on every fixed $\rho$-interval on which $F_*$ is defined, uniformly for $k$ near $\kappa$. Along this scaled orbit
$$
A=h^{-2}e^{2G_h}\rho^4.
$$
In this paragraph, $s=\log\rho$.

Choose $s_0$ large enough that the stationary orbit $U_*(s)$ is close to the cone for $s\geq s_0$. At $s_0$, the genuine orbit differs from $U_*$ by $O(h^{-2})$, and $A$ is comparable to $h^{-2}$. These statements remain true after subtracting $\mathcal Q_k$, uniformly for
$|k-\kappa|\leq Kh^{-1/2}$, with fixed $K$.

The variational equation along $U_*$ has a propagator bounded by
$Ce^{-\frac{\kappa+2}{2}(s-t)}$ for $s\geq t\geq s_0$. Indeed, the constant matrix $\mathsf J$ has this bound, and
$D\mathcal R_0(U_*(t))=O(e^{-\frac{\kappa+2}{2} t})$ is integrable. Variation of constants and Gronwall's inequality introduce only a bounded factor.
Consequently the difference from the stationary orbit satisfies
\begin{align*}
|V(s)-U_*(s)|\leq{}&Ch^{-2}e^{-\frac{\kappa+2}{2}(s-s_0)}\\
&+C\int_{s_0}^s e^{-\frac{\kappa+2}{2}(s-t)}
\big(A(t)|U_*(t)|+A(t)|V(t)-U_*(t)|+|V(t)-U_*(t)|^2\big)\,dt.
\end{align*}
As long as the orbit remains in the chosen neighborhood,
$A_s\geq (\kappa+2)A$. The preceding inequality closes the estimate
$$
|V(s)-U_*(s)|\leq C(h^{-2}+A(s))e^{-\frac{\kappa+2}{2} s}
$$
up to the first time $A=h^{-1}$. To check the bootstrap, use
$\int_{s_0}^s A(t)\,dt\leq CA(s)$. The term containing $A|U_*|$ then has the required size. The terms containing $A|V-U_*|$ and $|V-U_*|^2$ gain a small factor from $A\leq h^{-1}$ and the chosen neighborhood. Taking the bootstrap constant large and then $h$ large gives a strict improvement. The estimate also prevents an ordinary differential equation obstruction before this section; the positive lower bound for $A_s/A$ shows that the section is reached.

Integrating the first component, and using $\mathcal Q_k=O(|k-\kappa|A)$, gives
$$
G_h(\rho)-F_*(\rho)
=O(h^{-2})+O(|k-\kappa|A)+O(Ae^{-\frac{\kappa+2}{2} s}).
$$
Let $s_h$ be the first time $A=h^{-1}$. Substituting the logarithmic asymptotic of $F_*$ into $A=h^{-2}e^{2G_h}\rho^4$ yields
$$
s_h=\frac{\log h}{2(\kappa+2)}-\frac{c_{\rm in}}{(\kappa+2)}
+O(h^{-1/4}),\qquad
e^{-\frac{\kappa+2}{2} s_h}
=h^{-1/4}e^{c_{\rm in}/2}(1+O(h^{-1/4})).
$$
In particular, the matching phase comes from the fixed number $(\kappa+2)$, not from $k+2$.

\medskip\noindent\textbf{Propagation to a fixed outer section.}
In the variable $A$, the equation for $V$ is
$$
V_A=\left(\frac{\mathsf J}{2(\kappa+2)A}+\mathcal B_k(A)\right)V
+O(|V|^2/A),
$$
where $\mathcal B_k$ is bounded and analytic near zero, smoothly in $k$.
The linear equation has a fundamental matrix with columns
$$
\Phi_{k,0}(A)=A^{-1/2}(v_0+O(A)),\qquad
\Phi_{k,\pm}(A)=A^{-1/4\pm i\frac{\Omega}{2(\kappa+2)}}(v_\pm+O(A)).
$$
Here $v_0$ is a nonzero real eigenvector of $\mathsf J$ for $-(\kappa+2)$, $v_+$ is the eigenvector used in the stationary asymptotic expansion, and $v_-=\overline{v_+}$.
Indeed, the exponents have no nonzero integer differences. Substitution of a power series gives recursion matrices
$(\lambda+n)I-\mathsf J/(2(\kappa+2))$, $n\geq1$, whose inverses are $O(n^{-1})$.
A convergent majorant proves convergence, uniformly for $k$ in a smaller neighborhood. Extension over a fixed positive interval gives
$$
\|\Phi_k(A)\Phi_k(B)^{-1}\|\leq C(B/A)^{1/4},
\qquad 0<B\leq A\leq1.
$$

At $B=h^{-1}$, the stationary oscillatory expansion, the preceding comparison with $U_*$, and the formula for $s_h$ give a fixed nonzero complex number $c_0$ such that
$$
V(B)=h^{-1/2}\operatorname{Re}
\left(c_0e^{i(\Omega/(\kappa+2))\log h}\Phi_{k,+}(B)\right)+R_h,
\qquad |R_h|\leq Ch^{-1/2}\log h.
$$
The leading term has size $h^{-1/4}$. The error includes the real resonance term in the inner expansion and the $O(h^{-1/4})$ uncertainty in $s_h$.
Variation of constants then closes
$$
|V(A)|\leq Ch^{-1/2}A^{-1/4},\qquad h^{-1}\leq A\leq1.
$$
For example, its nonlinear contribution at $A=1$ is bounded by
$$
C\int_{h^{-1}}^1 A^{1/4}\frac{|V(A)|^2}{A}\,dA
\leq Ch^{-1}\int_{h^{-1}}^1 A^{-5/4}\,dA
\leq Ch^{-3/4}.
$$
At a general $A$, the same computation gives a relative gain at most
$Ch^{-1/4}$, which proves the bootstrap. The transported matching error is
$O(h^{-3/4}\log h)$. Replacing $\Phi_{k,+}(1)$ by
$\Phi_{\kappa,+}(1)$ costs $O(h^{-1})$ in the range
$|k-\kappa|\leq Kh^{-1/2}$. We have proved
\begin{align}
U(1;h,k)={}&\mathcal Q_k(1)
+h^{-1/2}\operatorname{Re}
\left(c_0e^{i(\Omega/(\kappa+2))\log h}\Phi_{\kappa,+}(1)\right)\notag\\
&+O(h^{-3/4}\log h), \label{eq:profile-passage}
\end{align}
uniformly in that range. Only a uniform continuous error estimate is needed below.

\medskip\noindent\textbf{Solving the two tail conditions.}
Let $\Psi_k:\mathbb R^3\to\mathbb R^2$ be a smooth local defining map, of rank two, for the tail curve at $A=1$ in the $U$ coordinates. It vanishes at $U=0$ when $k=\kappa$. Equation \eqref{eq:profile-passage} gives
$$
\Psi_k(U(1;h,k))
=h^{-1/2}\mathsf M
\begin{pmatrix}\cos(\frac{\Omega}{\kappa+2}\log h+\vartheta)\\\sin(\frac{\Omega}{\kappa+2}\log h+\vartheta)\end{pmatrix}
+(k-\kappa)B+O(h^{-3/4}\log h).
$$
Here $\mathsf M$ is a real $2$ by $2$ matrix, and $\vartheta$ is fixed.
The first linear connection assertion shows that $\mathsf M$ is invertible:
a nonzero vector in its kernel would give a pure oscillatory Frobenius solution lying on the tail line, contrary to the nonzero real Frobenius coefficient of $D$.

The vector $B$ includes the movement of the tail curve:
$$
B=\left.\partial_k\Psi_k(0)\right|_{k=\kappa}
+D_U\Psi_\kappa(0)
\left.\partial_k\mathcal Q_k(1)\right|_{k=\kappa}.
$$
It is nonzero. Otherwise the derivative of the curve $\mathcal Q_k$ would satisfy the linearized moving-tail condition. Its $p$ component $g$ solves
$(P(\theta)+x)g=x$ and is analytic and vanishes at zero. Differentiating the tail contraction gives
$g-1,g_q,g_{qq}=O(q^{-3})$. Since $x$ is a positive constant times $q^3$, this implies
$\theta g=O(q^{-2})$ and $\theta^2g=O(q^{-1})$ at infinity. Thus
$g,\theta g,\theta^2g$ are bounded at both endpoints, contrary to the adjoint-integral obstruction just proved. Therefore $B\ne0$.

Choose a nonzero real linear functional $\ell$ with $\ell(B)=0$.
Invertibility of $\mathsf M$ makes the function
$$\theta\longmapsto\ell(\mathsf M(\cos\theta,\sin\theta)^{\mathsf T})$$
a nonzero sinusoid.
At a simple zero $\theta_*$ there is a unique real $\nu_*$ such that
$$
\mathsf M\begin{pmatrix}\cos\theta_*\\\sin\theta_*\end{pmatrix}+\nu_*B=0.
$$
The map
$$
(\theta,\nu)\longmapsto\mathsf M\begin{pmatrix}\cos\theta\\\sin\theta\end{pmatrix}+\nu B
$$
has an invertible derivative there: applying $\ell$ to its $\theta$ derivative gives the nonzero derivative of the sinusoid, while its $\nu$ derivative is $B$. Choose a small closed disk about $(\theta_*,\nu_*)$ on whose boundary the displayed map does not vanish. Its degree is $1$ or $-1$; we use homotopy invariance and the existence property of the Brouwer degree as in~\cite[Chapter~1]{Deimling85}.

For large integers $j$, set
$$
h_j(\theta)=\exp\!\left(\frac{(\kappa+2)}{\Omega}
(\theta+2\pi j-\vartheta)\right),\qquad
k_j(\theta,\nu)=\kappa+h_j(\theta)^{-1/2}\nu.
$$
Choose $K$ larger than the maximum of $|\nu|$ on the disk. The maps
$$
h_j(\theta)^{1/2}\Psi_{k_j(\theta,\nu)}(U(1;h_j(\theta),k_j(\theta,\nu)))
$$
converge uniformly on that disk to the displayed map, since their error is
$O(h_j^{-1/4}\log h_j)$. They have the same nonzero degree for all large $j$, and hence have zeros. These zeros give $h_j\to\infty$ and
$k_j=\kappa+O(h_j^{-1/2})$. In particular $k_j>0$.
This argument does not require differentiating the matching error.

For the selected parameters, the genuine regular-origin orbit joins the constructed tail. Near the origin, in the scaled variable $\rho$,
$p/\rho^2\to1/6$. On every remaining fixed inner interval, $p>0$ follows from continuous dependence on the strictly positive stationary orbit. Past a sufficiently large fixed inner radius, the matching estimates keep $p$ close to $\kappa$; the tail keeps it close to $k_j>0$. Thus $p>0$ throughout. Since $A_s=2(p+2)A>0$, the reconstruction
$$
\frac{d\log r}{d\log A}=\frac1{2(p+2)},\qquad
F=\frac12\log A-2\log r
$$
gives a smooth profile for all $r>0$. Its smooth origin and the values
$F(0)=0$, $\Delta F(0)=h_j$ are retained, since the orbit used in the matching was the genuine regular-origin orbit.

The same matching gives $z>0$ everywhere for all sufficiently large selected $h_j$.
Near the regular origin this follows from $\Delta F(0)=h_j>0$.
On each remaining fixed interval in $\rho=\sqrt{h_j}r$, convergence
in the regular-origin construction gives positivity from the stationary orbit.
Choose the fixed inner section where that orbit is close to
$(\kappa,4\kappa,0)$. Up to $A=h_j^{-1}$ the comparison with $U_*$
keeps $z$ positive. On $h_j^{-1}\le A\le1$ the estimate
$|V(A)|\le Ch_j^{-1/2}A^{-1/4}\le Ch_j^{-1/4}$ and the smallness of
$\mathcal Q_k$ keep $z$ close to $4\kappa$.
The tail from $A=1$ to a fixed large $q_0$ is a small perturbation
of the exact cone, and the contracting tail satisfies
$z=4p+\tfrac23(p+2)qp_q=4k+o(1)$ with small uniform error.
These regions cover the orbit. Hence $\Delta F_j=z/r^2>0$.
Smoothness at zero, positivity on compact annuli and
$r^2\Delta F_j\to4k_j>0$ give the lower bound in
Lemma~\ref{lem:profile}.

\medskip\noindent\textbf{Differentiated tail expansions.}
We finish by proving the stated expansion to every finite order.
The tail system $X_q=B_kX+\mathcal F(q,k,X)$ is analytic in $q^{-1}$ and $X$ near $(0,0)$; its constant matrix $B_k$ is invertible.
For every $N\geq3$, substitution determines a formal vector sum
$$
X_N(q)=\sum_{n=3}^N d_nq^{-n}
$$
whose residual is $O(q^{-N-1})$. At each order, the coefficient of the new $d_n$ is $B_k$, and all remaining terms are already determined. The first two rows enforce that the second and third components are derivatives of the first.

Here the formal expansion is an asymptotic expansion of the actual tail. Subtract $X_N$ from the known solution $X=O(q^{-3})$.
The remainder solves a linear equation with coefficient
$$
B_k+\int_0^1D_X\mathcal F(q,k,X_N+t(X-X_N))\,dt,
$$
and source $O(q^{-N-1})$. The added matrix is $O(q^{-1})$. Its unstable terminal term at infinity vanishes because the known remainder is already $O(q^{-3})$. Write its stable and unstable variation-of-constants formula with the same stable initial component at a sufficiently large tail section. In the norm
$\sup q^{N+1}|X-X_N|$, the exponential convolution kernels are bounded, and the $O(q^{-1})$ coefficient is absorbed when that section is large. This constructs a remainder with the improved bound. It also belongs to the original $q^{-3}$ space, in which the same linear integral equation is unique. It must therefore equal the actual remainder. We conclude that
$$
X-X_N=O(q^{-N-1}).
$$
Differentiating the differential equation, after taking as many additional terms as needed, gives the expansion with every prescribed finite number of derivatives. This deduction uses the equation and additional terms, rather than differentiating an undifferentiated big-$O$ estimate.

In the scalar equation \eqref{eq:profile-tail-equation}, the differential terms lower a monomial by a total of three powers of $q$, the explicit source begins with $q^{-3}$, and the coefficient of $p-k$ at infinity is nonzero. The recursion therefore makes every scalar coefficient of $p-k$ outside the powers $q^{-3\ell}$ vanish. Equivalently, $p-k$ has an asymptotic series in $A^{-1}$.
Since $p-k=O(A^{-1})$ and
$A_s=2(p+2)A$, integration first gives
$$
F=k\log r+c_\infty+O(r^{-2(k+2)}).
$$
To pass to all orders, put $u=F-k\log r-c_\infty$ and use
$$
A=e^{2c_\infty}r^{2(k+2)}e^{2u},\qquad
\frac{du}{d\log r}=p-k.
$$
Substitution of the series for $p-k$ determines a series for $u$ in
$r^{-2(k+2)}$. Its recursion divides the coefficient of the $\ell$th term by
$-2(k+2)\ell$, which never vanishes. The already proved differentiated remainder estimates in $q$, and the identity
$d\log q/d\log r=2(p+2)/3$, give the corresponding remainder and all finite $r$ derivatives. This proves the asserted expansion.

Finally, $2a+4b=1$ and $a=bk$. Under
$f(t,x)=a\log(T-t)+F(|x|/(T-t)^b)$, both terms on the right-hand side of the evolution equation acquire the factor
$(T-t)^{-2a-4b}=(T-t)^{-1}$, while the left-hand side is
$(T-t)^{-1}(-a+brF_r)$. The profile equation therefore verifies the evolution equation exactly.
On every fixed physical annulus away from the origin, the tail
expansion becomes a series in integer powers of $T-t$ because
$2(k+2)b=1$. The differentiated expansions justify every finite
spatial and time derivative there, as required by the cutoff in
Section~\ref{sec:localization}.
\end{proof}

\section{The linearized operator}\label{sec:linear}

Fix a profile $F$ from Lemma~\ref{lem:profile}, and write
$g_{\mathrm E}=\sum_{j=1}^6(dy^j)^2$ for the Euclidean metric.
The auxiliary metric in this appendix is $g_F=2\gamma_F$; its
constant factor is used only for the analytic norms. Use the complete
metric $g_F=e^F g_{\mathrm E}$ on $\mathbb R^6$, complexifying the
function spaces for the spectral argument. The $L^2$ inner product
is linear in its first argument. At infinity set
$$
\rho=\frac{2e^{c_\infty/2}}{(k+2)}r^{(k+2)/2},\qquad P=e^{-F}\Delta.
$$
Extend $\rho$ smoothly through a compact set, and write
$\langle\rho\rangle=(1+\rho^2)^{1/2}$.
All adjoints below are taken with respect to
$d\mu_{g_F}=e^{3F}dy$. In particular,
$$
P^*u=e^{-3F}\Delta(e^{2F}u).
$$
The expansion in Lemma~\ref{lem:profile} gives, with derivatives
of every fixed order,
$$
g_F=(1+O(\rho^{-4}))\left(d\rho^2+\frac{(k+2)^2}{4}\rho^2g_{S^5}\right),\qquad br\partial_r=\frac\rho4\partial_\rho.
$$
Let $\nabla^F$ be the Chern connection of $g_F$; its constant
factor relative to $\gamma_F$ does not change the connection.
All higher covariant derivatives in this appendix use $\nabla^F$.
The scalar operator $\Delta_{g_F}$ below is the real Laplace--Beltrami
operator. The Chern and Levi--Civita connections differ by a tensor
linear in $dF$. The conformal connection formulas and the profile
expansion bound this tensor by $C\langle\rho\rangle^{-1}$ and its
$j$th covariant derivative by $C_j\langle\rho\rangle^{-1-j}$.
Consequently an $\ell$th derivative in either connection differs
from the other by a sum of lower derivatives of orders $j<\ell$
with coefficients $O(\langle\rho\rangle^{j-\ell})$.
Multiplication by $\langle\rho\rangle^{(m-\ell)/3}$ bounds each
such term by the weight $\langle\rho\rangle^{(m-j)/3}$.
Thus the weighted norms used below are equivalent in the two
connections. This justifies the real normal-coordinate estimates
while keeping all higher tensor norms in the Chern connection.

In particular, $g_F$ has bounded geometry, $P=\Delta_{g_F}-2\langle\nabla^{F} F,\nabla^{F}\cdot\rangle$, and $PF=O(\rho^{-2})$. Fix a sufficiently small $\alpha>0$, and choose a smooth function $\Phi$ equal to $\alpha\rho^{4/3}$ outside a compact set. For nonnegative integers $j$ define
$$
\|u\|_{\mathcal G_j}=\left(\sum_{\ell=0}^j\int_{\mathbb R^6}e^{-2\Phi}\langle\rho\rangle^{2(j-\ell)/3}|\nabla^{F,\ell}u|^2\,d\mu_{g_F}\right)^{1/2}.
$$
The spaces $\mathcal G_j$ are the completions of compactly supported smooth functions in these norms. Constants in this section may depend on $F$, $\Phi$, and the displayed derivative orders, but not on the length of a time interval.

\begin{lemma}\label{lem:linear}
The operator $\mathcal A=\mathcal N'_F-br\partial_r$ has a closed realization on $\mathcal G_0$ which generates an analytic semigroup and has compact resolvent. For sufficiently large fixed $\lambda_0>0$ and every nonnegative integer $j$,
$$
D((\lambda_0-\mathcal A)^j)=\mathcal G_{4j},\qquad (\mathcal G_{12},\mathcal G_{16})_{1/2,2}=\mathcal G_{14},
$$
with equivalent norms. Choose $\delta>0$ such that $\{\operatorname{Re}z=-\delta\}$ misses the spectrum. Let $\Pi_+$ be the spectral projection onto the eigenvalues to its right, and put $\Pi_-=I-\Pi_+$. Then $\Pi_+$ has finite rank.

For $I_\theta=[\max(\tau_0,\theta-1),\theta]$, define
\begin{align*}
\|u\|_{\mathfrak X_\delta[\tau_0,S]}&=\sup_{\tau_0\leq\theta\leq S}e^{\delta\theta}\left(\|u\|_{L^\infty(I_\theta;\mathcal G_{14})}+\|u\|_{L^2(I_\theta;\mathcal G_{16})}+\|u_\tau\|_{L^2(I_\theta;\mathcal G_{12})}\right),\\
\|h\|_{\mathfrak Y_\delta[\tau_0,S]}&=\sup_{\tau_0\leq\theta\leq S}e^{\delta\theta}\|h\|_{L^2(I_\theta;\mathcal G_{12})}.
\end{align*}
The initial trace is included when $I_\theta$ is a singleton. Solutions of $u_\tau=\mathcal A u+\Pi_-h$ with $u(\tau_0)\in\operatorname{ran}\Pi_-$ satisfy
$$
\|u\|_{\mathfrak X_\delta[\tau_0,S]}\leq C\left(e^{\delta\tau_0}\|u(\tau_0)\|_{\mathcal G_{14}}+\|h\|_{\mathfrak Y_\delta[\tau_0,S]}\right),
$$
where $C$ is independent of $S$. Finally, if $U=\operatorname{ran}\Pi_+$ and $\beta>0$, there is a linear map $\mathcal J:U\longrightarrow C_c^\infty(\mathbb R^6)$ such that $\Pi_+\mathcal J=I$ and $\|\Pi_-\mathcal J\|_{U\to\mathcal G_{14}}<\beta$.
\end{lemma}

\begin{proof}
Functions added to differential operators are understood as
multiplication operators; thus $(P+PF)u=Pu+(PF)u$.
The same convention applies to $9(PF)^2-3\mathcal N(F)$.
The linearized operator admits the exact factorization
$$
\mathcal N'_F=-18(P+PF)^*(P+PF)+9(PF)^2-3\mathcal N(F).
$$
To verify it without assigning a finite energy to $F$, differentiate the relative functional
$$
9\int_{\mathbb R^6}\bigl(e^w(PF+Pw)^2-(PF)^2\bigr)\,d\mu_{g_F}
$$
for compactly supported variations. Its Hessian at zero is
$$
18\langle(P+PF)u,(P+PF)v\rangle-9\langle(PF)^2u,v\rangle.
$$
Differentiating the factor $e^{-3w}$ in the negative gradient contributes $-3\mathcal N(F)u$, which gives the factorization. The multiplication coefficient $9(PF)^2-3\mathcal N(F)$ and all its derivatives are bounded; it is $O(\rho^{-4})$ at infinity because $\mathcal N(F)=b(rF_r-k)$.

Conjugate $-\mathcal A$ by $u\mapsto e^{-\Phi}u$ onto $L^2(g_F)$. Its form on compactly supported smooth functions is
\begin{align*}
\mathfrak a(u,v)={}&18\left\langle e^{-\Phi}(P+PF)e^\Phi u,e^\Phi(P+PF)e^{-\Phi}v\right\rangle+\langle br\partial_r u,v\rangle\\
&+\left\langle\bigl(br\partial_r\Phi-9(PF)^2+3\mathcal N(F)\bigr)u,v\right\rangle.
\end{align*}
Its form domain is $H^2(g_F)\cap\{\langle\rho\rangle^{2/3}u\in L^2(g_F)\}$. To prove coercivity, integration by parts gives
$$
\int\langle\rho\rangle^{2/3}|\nabla^{F} u|^2\,d\mu_{g_F}\leq C\|\Delta_{g_F}u\|_2\|\langle\rho\rangle^{2/3}u\|_2+C\|u\|_2^2.
$$
The two conjugated second order operators in the form are
$$
P+PF+|\nabla^{F}\Phi|^2\ \pm\ \bigl(2\langle\nabla^{F}\Phi,\nabla^{F}\cdot\rangle+P\Phi\bigr).
$$
The real part of their pairing is the squared norm of the common part minus the squared norm of the signed part. Since $\nabla^{F}\Phi=O(\alpha\langle\rho\rangle^{1/3})$, the preceding interpolation estimate and Young's inequality bound the latter loss by
$$
\varepsilon\|\Delta_{g_F}u\|_2^2+C_\varepsilon\alpha^4\|\langle\rho\rangle^{2/3}u\|_2^2+C_\varepsilon\|u\|_2^2.
$$
Also $br\partial_r\Phi=\alpha\rho^{4/3}/3$ at infinity, while the real part of $\langle br\partial_r u,u\rangle$ is bounded by $C\|u\|_2^2$ because the drift has bounded divergence. Taking $\alpha$ small gives
$$
\operatorname{Re}\mathfrak a(u,u)+C\|u\|_2^2\geq c\|\Delta_{g_F}u\|_2^2+\frac\alpha6\|\langle\rho\rangle^{2/3}u\|_2^2.
$$
The form is bounded on its domain. For the drift, split its coefficient $O(\rho)$ as $O(\langle\rho\rangle^{1/3})\langle\rho\rangle^{2/3}$ and use the same interpolation inequality. This also bounds the imaginary part of the form by a constant times its shifted real part. The representation theorem for closed sectorial forms~\cite[Chapter~VI, Theorem~2.1]{Kato95} therefore defines a sectorial operator with angle smaller than $\pi/2$. The rotated Lax--Milgram estimate supplies its resolvents, and the resolvent contour construction gives the analytic semigroup; see~\cite[Chapter~2]{Lunardi95}. The form-domain inclusion in $L^2$ is compact: Rellich's theorem applies on fixed balls, and the squared $L^2$ tail is at most $R^{-4/3}\|\langle\rho\rangle^{2/3}u\|_2^2$. This proves compactness of the resolvent.

We next identify its operator domain; form coercivity alone would not suffice for the nonlinear argument. For the conjugated operator the required estimate is
$$
\sum_{j=0}^4\|\langle\rho\rangle^{(4-j)/3}\nabla^{F,j} u\|_2\leq C\bigl(\|e^{-\Phi}\mathcal A(e^\Phi u)\|_2+\|u\|_2\bigr).
$$
Here is a localization proof with control of the domain at infinity. Let $\Lambda\ge1$ be a constant to be chosen. At a point with
$\rho=R\gg1$, use normal coordinates on a ball of radius
$\Lambda R^{-1/3}$ and freeze the coefficients. Writing $\eta=\nabla^{F}\Phi$ at its center, the terms of fourth order in this scale have symbol
$$
18\bigl(|\xi|^2-|\eta|^2-2i\eta\cdot\xi\bigr)^2+i\langle br\partial_r,\xi\rangle+\frac\alpha3\rho^{4/3}.
$$
The drift term is purely imaginary. The real part is bounded below by $c|\xi|^4+c\alpha R^{4/3}$: its negative mixed term is at most $108|\eta|^2|\xi|^2$, which is absorbed by half the $|\xi|^4$ term and $C\alpha^4R^{4/3}$. Fourier transformation consequently controls $R^{(4-j)/3}D^ju$ for $0\leq j\leq4$. On these balls the coefficient variation has relative graph bound $O(\Lambda R^{-4/3})+o(1)$. A partition with bounded overlap has derivatives of order $j$ bounded by $C_j\Lambda^{-j}R^{j/3}$, so its commutators, including the drift commutator, have relative graph bound $O(\Lambda^{-1})$. Choose $\Lambda$ large and then the starting radius large. Absorption, followed by ordinary local elliptic estimates on the remaining compact set, proves the displayed graph estimate for compactly supported functions.

Let now $u$ belong to the operator domain defined by the form. It is locally in $H^4$. For fixed $N>4/3$, put $\vartheta_R=(1+\langle\rho\rangle/R)^{-N}$. Interior elliptic estimates on balls of radius $\langle\rho\rangle^{-1/3}$ first show that $\vartheta_Ru$ has all the weighted derivatives in the graph estimate: the only additional local error is $\langle\rho\rangle^{4/3}u$, and $\vartheta_R\langle\rho\rangle^{4/3}$ is bounded for fixed $R$. Thus the estimate is applicable before any uniform bound is asserted. Expressing the commutator with $\vartheta_R$ in derivatives of $\vartheta_Ru$ gives a relative graph bound tending to zero as $R\to\infty$. For its leading term the ratio is
$$
\frac{C}{(R+\langle\rho\rangle)\langle\rho\rangle^{1/3}},
$$
and the drift gives $C\rho/((R+\langle\rho\rangle)\langle\rho\rangle^{4/3})$; all lower terms are smaller. Absorb these commutators and use Fatou's lemma to obtain the graph estimate for $u$. Conversely, smooth compactly supported approximation in this graph norm and integration by parts identify every such weighted Sobolev function with an element of the form realization. This proves domain equality. Undoing the conjugation preserves these norms, since each derivative of $e^\Phi$ supplies only the corresponding factor $O(\langle\rho\rangle^{1/3})$.

For completeness the argument also identifies the higher domains. The frozen symbol of a shifted $j$th power is the $j$th power of the preceding symbol and has modulus bounded below by $c_j(|\xi|^4+R^{4/3})^j$. Its coefficient of derivative order $\ell$ has scale at most $R^{(4j-\ell)/3}$, and coefficient derivatives have the same small normalized variations as above. The partition argument therefore gives all derivatives in $\mathcal G_{4j}$. For an element of the power domain, local elliptic regularity gives $H^{4j}_{\rm loc}$; the same soft cutoff with $N>4j/3$ makes its initial localized graph norm finite. The commutator divided by the relevant graph weight again tends uniformly to zero, so absorption and exhaustion prove the domain inclusion. The reverse inclusion follows by applying the differential operator successively, using the coefficient bounds and compactly supported approximation. Hence $D((\lambda_0-\mathcal A)^j)=\mathcal G_{4j}$.

The interpolation assertion can be checked in the same adaptive cover. On a ball of radius $\langle\rho\rangle^{-1/3}$, rescale to a unit ball. The $\mathcal G_j$ norm there is, up to a common volume factor, $\langle\rho\rangle^{j/3}e^{-\Phi}$ times the ordinary $H^j$ norm; $\Phi$ varies by a bounded amount on the ball. The localization and reconstruction maps are bounded for each of the finitely many relevant orders. Real interpolation thus reduces to the Fourier identity $(H^{12},H^{16})_{1/2,2}=H^{14}$, together with interpolation of the displayed scalar weights. This proves the claimed trace space without a square-root-domain assertion for a non-self-adjoint operator.

The sectorial resolvent bound and compactness imply that only finitely many eigenvalues lie in $\{\operatorname{Re}z\geq-\delta\}$. The resolvent contour on the complement of $\Pi_+$ yields
$$
\|e^{s\mathcal A}\Pi_-\|_{\mathcal G_j\to\mathcal G_j}\leq Ce^{-(\delta+\varepsilon)s},\qquad s\geq0,
$$
for some $\varepsilon>0$ and each of the graph and trace spaces used here. To obtain maximal regularity on the stable subspace of $\mathcal G_{12}$, shift by $\delta$ and take a time Fourier transform for zero initial data. The resolvent estimates bound both multipliers
$$
i\xi(i\xi-\mathcal A-\delta)^{-1}\Pi_-,\qquad (\mathcal A+\delta)(i\xi-\mathcal A-\delta)^{-1}\Pi_-.
$$
Plancherel gives the time derivative and $\mathcal G_{16}$ bounds,
and the one-dimensional interpolation trace inequality gives the
$\mathcal G_{14}$ trace. For a nonzero initial trace
$u_0\in\mathcal G_{14}$, analyticity and the domain estimate give,
for $0<s\le1$,
$$
\|\mathcal Ae^{s\mathcal A}\Pi_-u_0\|_{\mathcal G_{12}}
\le\frac Cs\inf_{u_0=u+v}
\left(\|u\|_{\mathcal G_{12}}+s\|v\|_{\mathcal G_{16}}\right).
$$
The infimum is over $u\in\mathcal G_{12}$ and
$v\in\mathcal G_{16}$. Squaring and integrating in $s$ gives the
$L^2\mathcal G_{16}$ and time-derivative bounds from the defining
real interpolation norm of
$(\mathcal G_{12},\mathcal G_{16})_{1/2,2}$.
Boundedness of the semigroup on the two endpoint spaces, followed
by interpolation, gives its $\mathcal G_{14}$ trace bound.
Exponential stability controls $s\ge1$.

For the stated uniformly local estimate, divide the forcing into
consecutive unit intervals. Its contribution on its own interval
has the preceding estimate; at every later interval its endpoint
trace is propagated by $e^{s\mathcal A}\Pi_-$. After multiplication
by $e^{\delta\theta}$ these contributions are bounded by the
summable sequence $Ce^{-\varepsilon n}$. Arbitrary backward unit
intervals meet at most two of these blocks. This proves the
estimate with a constant independent of $S$.

Finally, every generalized eigenfunction belongs to all power domains, hence is smooth. If $\chi_R$ is a smooth compactly supported cutoff tending to one, then $\chi_Re\to e$ in $\mathcal G_{14}$ for every $e\in U$. For large $R$, the map $\Pi_+\chi_R:U\to U$ is invertible. Taking $\mathcal J=\chi_R(\Pi_+\chi_R|_U)^{-1}$ gives $\Pi_+\mathcal J=I$, and $\Pi_-\mathcal J\to0$ in the asserted operator norm. The spectral set is invariant under conjugation, so the projections commute with complex conjugation. Taking real cutoffs restricts this construction to the real invariant subspace and gives real-valued initial perturbations. This completes the proof.
\end{proof}

\begin{lemma}\label{lem:tame}
For real smooth $w$, set $Q(w)=\mathcal N(F+w)-\mathcal N(F)-\mathcal N'_Fw$. If $\sum_{j=0}^8\|\nabla^{F,j}w\|_\infty\leq\varepsilon\leq1$, then
$$
\|Q(w)\|_{\mathcal G_{12}}\leq C\varepsilon\|w\|_{\mathcal G_{16}},\qquad \|\Pi_+Q(w)\|_U\leq C\varepsilon\|w\|_{\mathcal G_{14}}.
$$
The hypotheses include the displayed unweighted bounds; they are not a consequence of the exponentially weighted norm alone.
\end{lemma}

\begin{proof}
In the fixed metric $g_F$ the exact nonlinear expression is
$$
\mathcal N(F+w)=-18e^{-3w}P^*\bigl(e^w(PF+Pw)\bigr)-9e^{-2w}(PF+Pw)^2.
$$
All background coefficients and their covariant derivatives are bounded. After subtracting the constant and linear terms, every product contains at least two factors involving $w$. After $\ell\leq12$ further derivatives, their total derivative order is at most $\ell+4\leq16$. Put a factor of highest order $j$ in weighted $L^2$. Every other factor has order at most eight, since two factors of order at least nine would exceed total order sixteen. They are bounded by the hypothesis, while
$$
\langle\rho\rangle^{(12-\ell)/3}\leq\langle\rho\rangle^{(16-j)/3},\qquad j\leq\ell+4.
$$
Leibniz's rule and the chain rule prove the first estimate. A coefficient such as $e^{-2w}-1$ without derivatives is $O(\varepsilon)$; its differentiated terms obey the same count. At order zero the argument already gives $\|Q(w)\|_{\mathcal G_0}\leq C\varepsilon\|w\|_{\mathcal G_4}$. Boundedness of the finite-rank projection on $\mathcal G_0$ and the inclusion $\mathcal G_{14}\subset\mathcal G_4$ prove the second estimate. No pointwise $\mathcal G_{16}$ bound and no algebra property for the exponential weight is used.
\end{proof}

\vspace{0.3cm}

{\small{\bf{Acknowledgments:}} {Part of this work was done during Ma's visit to Tianyuan mathematical center in southwest China in June of 2025.
}}

\vspace{0.3cm}

%\textbf{Funding} No funding.

\textbf{Data Availability Statement}
Data sharing is not applicable to this article because no datasets were generated or analyzed during the current study.

\textbf{Compliance with ethical standards} %Not applicable.

\textbf{Conflict of interest} The authors declare that they have no conflict of interest.

\textbf{Code availability}
Not applicable.

\textbf{AI Assistance Statement}
AI was used to assist with this research. The authors take full responsibility for the mathematical content and conclusions.


\begin{thebibliography}{CMS26}
\setlength{\itemsep}{2pt}

\bibitem[ACS17]{ACS17}
D. Angella, S. Calamai and C. Spotti,
On the Chern--Yamabe problem,
\emph{Mathematical Research Letters} \textbf{24} (2017), no.~3,
645--677.

\bibitem[BV20]{BV20}
L. Bedulli and L. Vezzoni,
A scalar Calabi-type flow in Hermitian geometry: short-time existence
and stability,
\emph{Annali della Scuola Normale Superiore di Pisa, Classe di Scienze}
(5) \textbf{20} (2020), no.~2, 657--676.

\bibitem[CH08]{CH08}
X. X. Chen and W. Y. He,
On the Calabi flow,
\emph{American Journal of Mathematics} \textbf{130} (2008), no.~2,
539--570.

\bibitem[CMS26]{CMS26}
H. Chen, L. Ma and L. Sheng,
\emph{Continuation and convergence for the higher dimensional
conformal Chern--Calabi flow},
\href{https://arxiv.org/abs/2610.31270v1}{arXiv:2610.31270v1} (2026).

\bibitem[D85]{Deimling85}
K. Deimling,
\emph{Nonlinear Functional Analysis},
Springer-Verlag, Berlin--Heidelberg, 1985, Chapter~1.

\bibitem[DK10]{DK10}
H. Dong and D. Kim,
\emph{On the $L_p$-solvability of higher order parabolic and elliptic
systems with BMO coefficients},
\href{https://arxiv.org/abs/1003.0969v2}{arXiv:1003.0969v2} (2010).

\bibitem[DZ15]{DZ15}
H. Dong and H. Zhang,
Schauder estimates for higher-order parabolic systems with time
irregular coefficients,
\emph{Calculus of Variations and Partial Differential Equations}
\textbf{54} (2015), no.~1, 47--74.
Version used for Theorem~2.2(I):
\href{https://arxiv.org/abs/1307.5013v1}{arXiv:1307.5013v1}.

\bibitem[EMT11]{EMT11}
J. Enders, R. M\"uller and P. M. Topping,
On Type-I singularities in Ricci flow,
\emph{Communications in Analysis and Geometry} \textbf{19} (2011),
no.~5, 905--922.

\bibitem[K95]{Kato95}
T. Kato,
\emph{Perturbation Theory for Linear Operators}, second edition,
reprint, Classics in Mathematics, Springer-Verlag, Berlin, 1995.

\bibitem[K26]{Kawamura26}
M. Kawamura,
On the Calabi flow on quasi-K\"ahler manifolds,
\emph{Partial Differential Equations and Applications}
\textbf{7} (2026), article~1.

\bibitem[L95]{Lunardi95}
A. Lunardi,
\emph{Analytic Semigroups and Optimal Regularity in Parabolic Problems},
Progress in Nonlinear Differential Equations and Their Applications,
vol.~16, Birkh\"auser, Basel, 1995.

\bibitem[MM12]{MM12}
C. Mantegazza and L. Martinazzi,
A note on quasilinear parabolic equations on manifolds,
\emph{Annali della Scuola Normale Superiore di Pisa, Classe di Scienze}
(5) \textbf{11} (2012), no.~4, 857--874.

\bibitem[S22]{Shen22}
X. S. Shen,
A Chern--Calabi flow on Hermitian manifolds,
\emph{Journal of Geometric Analysis} \textbf{32} (2022), no.~4,
article~129.

\bibitem[SM26]{MS26}
L. Sheng and L. Ma,
\emph{On global Chern-Calabi flows on compact complex surfaces},
\href{https://arxiv.org/abs/2610.31068v1}{arXiv:2610.31068v1} (2026).




\bibitem[S13a]{Streets13long}
J. Streets,
The long time behavior of fourth order curvature flows,
\emph{Calculus of Variations and Partial Differential Equations}
\textbf{46} (2013), nos.~1--2, 39--54.

\bibitem[S13b]{Streets13collapse}
J. Streets,
Collapsing in the $L^2$ curvature flow,
\emph{Communications in Partial Differential Equations}
\textbf{38} (2013), no.~6, 985--1014.

\bibitem[Y25]{Yudowitz25}
L. Yudowitz,
Semi-continuity of the Morse index for Ricci shrinkers,
\emph{Journal of Geometric Analysis} \textbf{35} (2025), article~159.

\end{thebibliography}
\end{document}